\documentclass[12pt]{amsart}
\usepackage{amssymb}
\newtheorem{theorem}{Theorem}[section]
\newtheorem{lemma}[theorem]{Lemma}
\newtheorem{proposition}[theorem]{Proposition}
\newtheorem{corollary}[theorem]{Corollary}

\theoremstyle{definition}
\newtheorem{definition}[theorem]{Definition}

\theoremstyle{remark}
\newtheorem{remark}[theorem]{Remark}
\numberwithin{equation}{section}

\newcommand{\N}{\mathbb{N}}
\newcommand{\Z}{\mathbb{Z}}

\newcommand{\R}{\mathbb{R}}
\newcommand{\C}{\mathbb{C}}
\newcommand{\T}{\mathbb{T}}
\newcommand{\tM}{\widetilde{M}}
\newcommand{\tN}{\widetilde{N}}
\newcommand{\ta}{\widetilde{\alpha}}
\newcommand{\tr}{\widetilde{\rho}}
\newcommand{\tm}{\widetilde{\mu}}
\newcommand{\tE}{\widetilde{E}}
\newcommand{\tJ}{\widetilde{J}}
\newcommand{\tb}{\widetilde{\beta}}
\newcommand{\tP}{\widetilde{P}}
\newcommand{\tQ}{\widetilde{Q}}
\newcommand{\tep}{\widetilde{\epsilon}}
\newcommand{\eM}{\mathrm{End}(M)_0}
\newcommand{\eN}{\mathrm{End}(N)_0}
\newcommand{\eMm}{\mathrm{End}(M)_m}
\newcommand{\eNm}{\mathrm{End}(N)_m}
\newcommand{\eMC}{\mathrm{End}(M)_{CT}}
\newcommand{\eMt}{\mathrm{End}(M)_{tr}}
\newcommand{\sM}{\mathrm{Sect}(M)_0}

\newcommand{\sMm}{\mathrm{Sect}(M)_m}

\newcommand{\sMC}{\mathrm{Sect}(M)_{CT}}
\newcommand{\sMt}{\mathrm{Sect}(M)_{tr}}
\newcommand{\Di}{\mathrm{Diag}}
\newcommand\inpr[2]{\langle{#1,#2}\rangle}
\newcommand\id{\mathrm{id}}
\newcommand\br{\overline{\rho}}

\newcommand\bc{\overline{c}}
\newcommand\bR{\overline{R}}
\newcommand\bp{\overline{\pi}}
\newcommand\mo{\mathrm{mod}}
\newcommand\Mo{\mathrm{Mod}}
\newcommand\Ad{\mathrm{Ad}}
\newcommand\hM{\widehat{M}}
\newcommand\hr{\widehat{\rho}}
\newcommand\hK{\widehat{K}}
\newcommand\hE{\widehat{E}}
\newcommand\hb{\widehat{\beta}}
\newcommand\cR{\mathcal{R}}
\newcommand\cH{\mathcal{H}}
\newcommand\cK{\mathcal{K}}

\newcommand\cD{\mathcal{D}}
\newcommand\cE{\mathcal{E}}
\newcommand\cL{\mathcal{L}}
\begin{document}

\title{Canonical Extension of Endomorphisms of Type III Factors}

\author{Masaki Izumi}

\address{Department of Mathematics\\ Graduate School of Science\\
Kyoto University\\ Sakyo-ku, Kyoto 606-8502\\ Japan}
\email{izumi@kusm.kyoto-u.ac.jp}

\thanks{Work supported by Mathematical Sciences Research Institute}

\begin{abstract} We extend the notion of the canonical extension of
automorphisms of type III factors to the case of endomorphisms with
finite statistical dimensions.
Following the automorphism case, we introduce two notions for
endomorphisms of type III factors: modular endomorphisms and Connes-Takesaki
modules.
Several applications to compact groups of automorphisms and subfactors of
type III factors are given from the viewpoint of ergodic theory.
\end{abstract}

\maketitle

\section{Introduction}
The canonical extension of an automorphism of a type III factor gives
rise to a natural isomorphism from the automorphism group of the type III
factor into that of the crossed product by the modular automorphism group.
Although the notion itself had been known to specialists before,
(e.g. \cite{D}), the first systematic analysis of the canonical
extension was proposed and accomplished by U. Haagerup and E. St{\o}rmer.
In \cite{HS1}, \cite{HS2}, they introduced several classes of automorphisms
of von Neumann algebras and investigated their structure using
the canonical extension.

One of the purposes of the present notes is to generalize the notion
of the canonical extension to an endomorphism whose image has a finite
index \cite{J2}, \cite{PP}, \cite{K1}.
Our analysis is based on techniques of the common continuous decomposition
of inclusions and endomorphisms of type III factors, mainly developed by
Ph. Loi \cite{Loi}, T. Hamachi-H. Kosaki \cite{HK1}, \cite{HK2},
\cite{HK3}, R. Longo \cite{L2}, \cite{L3}, and H. Kosaki-R. Longo \cite{KL}.
In \cite{R}, J. E. Roberts introduced an action of the dual object of a
compact group on a von Neumann algebra, now called a Roberts action,
which is a functor from the category of the finite dimensional unitary
representations to that of the endomorphisms of the von Neumann algebra.
As several invariants of discrete group actions on type III factors can be
captured by the canonical extension \cite{ST1}, \cite{Se},
it is natural to expect that the canonical extension of endomorphisms
plays an important role in the analysis of compact group actions through
the dual Roberts actions.
Indeed, we obtain several new results on minimal actions of compact groups
on type III factors.
Study of minimal actions of compact groups was initiated by S. Popa and
A. Wassermann using subfactor techniques developed in \cite{P}.
Interested readers are recommended to consult \cite{PW} for related topics.

We introduce two new notions for endomorphisms corresponding to
extended modular automorphisms and the Connes-Takesaki module
in the automorphism case \cite{CT}: modular endomorphisms and
the Connes-Takesaki module of endomorphisms.
As in the case of extended modular automorphisms,
a modular endomorphism carries a unitary group valued cocycle
of the flow of weights, and it gives rise to a Roberts action
of the dual object of the essential range of the cocycle,
the notion developed by G. Mackey and R. Zimmer \cite{M3}, \cite{Z1}.
The dual actions of Roberts actions consisting of modular endomorphisms
provide a large variety of new examples of minimal actions of compact
groups on type III$_0$ factors.
Also,  we give a strong constraint for the possible types of
an algebra-fixed point algebra pair for a minimal action of a compact
semisimple Lie group using these two notions.

In \cite{Z1}, Zimmer introduced the notion of extensions of ergodic
transformation groups with relatively discrete spectrum,
which is a relative version of the notion of ergodic transformations
with pure point spectrum.
It turns out that a non-commutative generalization of this notion is
closely related to modular endomorphisms (Theorem 5. 11), which may be
expected because both objects are described in terms of cocycles of
ergodic transformation groups with values in compact groups.
In particular, we characterize a subfactor, not necessarily of
a finite index, coming from an ergodic extension of the flow of
weights in terms of a group-subgroup subfactor
(c.f. \cite{Ham1}, \cite{Ham2}, \cite{HK1}, \cite{HK2}, \cite{Su2}).

Our original purpose for introducing the canonical extension of
endomorphisms is to settle a problem left unsolved in \cite{I1}.
When a type III factor and its subfactor with a finite index have
the common flow of weights, their ``type II principal graph'' makes
sense via the common crossed product decomposition
\cite{Loi}, \cite{KL}.
In \cite{I1}, we characterized the situation where the type II
principal graph does not coincide with the principal graph of
the original subfactor for the case of III$_\lambda$, $\lambda\neq 0$
factors.
Namely, the two graphs do not coincide if and only if some power of the
canonical endomorphism for the inclusion contains a non-trivial
modular automorphism.
It would be tempting to conjecture that the same statement should be
true in the type III$_0$ case if the modular automorphism is replaced
with an extended modular automorphism.
However, it turns out that the right counterpart for our purpose
is a modular endomorphism (Theorem 3.7).
One can observe a similar phenomenon in the recent work of
T. Masuda \cite{Mas} (c.f. \cite{K3}).

Final part of this work was finished while the author stayed at MSRI,
and he would like to thank them for their hospitality.
He also would like to thank T. Hamachi and  H. Kosaki for stimulating
discussions.
Some of the main results in this paper were announced in \cite{I2}.

\section{Canonical Extension}
First, we briefly summarize notation used in this paper.
Our basic references are \cite{St} and \cite{SZ1} for Tomita-Takesaki
theory, \cite{K1} for index theory of type III factors,
\cite{ILP} for sector theory and infinite index inclusions, and
\cite{Z1} for cocycles of ergodic transformation groups.
Undefined terms and notations used in this paper should be found
in these references.
We always assume that von Neumann algebras have separable preduals,
Hilbert spaces are separable, and locally compact groups are
second countable.
Automorphisms and endomorphisms of von Neumann algebras are always
assumed to be $*$-preserving and unital.
For a von Neumann algebra $M$,
we denote by $\mathrm{Aut}(M)$, $\mathrm{Int}(M)$, and $\mathrm{End}(M)$
the sets of automorphisms, inner automorphisms, and endomorphisms
respectively.
We denote by $\mathrm{Out}(M)$ and $\mathrm{Sect}(M)$ the sets of
unitary equivalence classes of automorphisms and endomorphisms
respectively.
A member in $\mathrm{Sect}(M)$ is called a sector.
For $\rho, \sigma\in \mathrm{End}(M)$, the intertwiner space
from $\rho$ to $\sigma$ is defined by
$$(\rho,\sigma)=\{V\in M;\; V\rho(x)=\sigma(x)V\}.$$
When $M$ is a factor and $\rho$ is irreducible, that is,
$\rho(M)'\cap M=\C 1$,
$(\rho,\sigma)$ is a Hilbert space with the  inner product
$$\inpr {V_1}{V_2} 1=V_2^*V_1.$$
Every action of a topological group on a von Neumann algebra $M$
is assumed to be continuous with respect to the $u$-topology in
$\mathrm{Aut}(M)$.
For a group action $\alpha$ on $M$, $M^\alpha$ denotes the fixed point
subalgebra of $M$ under the action $\alpha$.
When $M$ acts on a Hilbert space $H$, we denote by $\pi_\alpha$
the representation of $M$ on $L^2(G,H)$ defined by
$$\pi_\alpha(x)\xi(g)=\alpha_{g^{-1}}(x)\xi(g),\quad \xi\in L^2(G,H).$$
We denote by $l(g)$ and $r(g)$ the left and right regular representations
of $G$ respectively.
For a weight $\varphi$ on $M$ or an operator valued weight from $M$ to
a subalgebra, we denote by $\mathfrak{m}_\varphi$ and
$\mathfrak{n}_{\varphi}$ the domain of $\varphi$ and the left ideal
corresponding to $\mathfrak{m}_\varphi$ as usual.

Let $\varphi$ be a faithful normal semifinite weight on $M$.
We denote by $\tM=M\rtimes_{\sigma^\varphi}\R$ the crossed
product of $M$ by the modular automorphism group
$\{\sigma^\varphi_t\}_{t\in \R}$, which is the von Neumann algebra
generated by $\pi_{\sigma^{\varphi}}(M)$ and the implementing
one-parameter unitary group
$\{\lambda^\varphi(t)\}_{t\in \R}$, where
$\lambda^\varphi(t)=1\otimes l(t)$.
We often omit $\pi_{\sigma^{\varphi}}$ when there is no possibility of
confusion.
We denote by $\theta$ and $\tau$ the dual action of $\sigma^\varphi$
and the natural trace constructed from the dual weight of $\varphi$ on
$\tM$ and the generator of $\{\lambda^\varphi(t)\}_{t\in \R}$.
Then, the triple $(\tM, \theta, \tau)$ does not depend on the choice of
$\varphi$ under the following identification
$$\lambda^\psi(t)=[D\psi:D\varphi]_t\lambda^\varphi(t),$$
where $\psi$ is another faithful normal weight on $M$ and
$\{[D\psi:D\varphi]_t\}_{t\in \R}$ is the Connes cocycle derivative.
As for the relationship between $\mathrm{Aut}(M)$ and $\mathrm{Aut}(\tM)$,
we have the following  \cite[Proposition 12.1]{HS2}:

\begin{proposition} Let $M$ be a von Neumann algebra and $\varphi$ be a
faithful normal semifinite weight on $M$.
Then, for every automorphism $\alpha\in \mathrm{Aut}(M)$,
there exists a unique automorphism $\ta\in \mathrm{Aut}(M)$
satisfying
$$\ta(x)=\alpha(x),\quad x\in M,$$
$$\ta(\lambda^\varphi(t))=[D\varphi\cdot\alpha^{-1}:D\varphi]_t
\lambda^\varphi(t).$$
Moreover, the map $\alpha\mapsto \ta$ is a homomorphism from
$\mathrm{Aut}(M)$
into $\mathrm{Aut}(\tM)$.
\end{proposition}

$\ta$ is called the canonical extension of $\alpha$.
The goal of this section is to introduce the canonical extension
of endomorphisms with finite statistical dimensions generalizing
the automorphism case.

Let $M$ be an infinite factor.
We denote by $\mathrm{End}(M)_0$ the set of endomorphisms of $M$
with finite statistical dimensions.
For $\rho\in \eM$, we denote by $E_\rho$, $\phi_\rho$, and $d(\rho)$,
the minimal conditional expectation from $M$ onto the image of $\rho$
\cite{H}, the standard left inverse of $\rho$, and the statistical
dimension of $\rho$, respectively:
that is, $\phi_\rho=\rho^{-1}\cdot E_\rho$ and
$d(\rho)=[M:\rho(M)]_0^{1/2}$, where
$[M:\rho(M)]_0=\mathrm{Index}\; E_\rho$ is the minimum index of
$M\supset \rho(M)$.
We denote by $\sM$ the unitary equivalence classes of the endomorphisms
in $\eM$.
$\mathrm{Sect}(M)$ has three natural operations,
forming a direct sum, composing two endomorphisms (regarded as a product),
and taking conjugation.
The statistical dimensions are additive, multiplicative,
and invariant under these three operations on $\sM$ respectively
\cite{H}, \cite{L3}, \cite{KW}, \cite{LR}.
Let $M \supset N$ be an inclusion of type III factors with a finite index.
Then, we regard $\tN$ as a subalgebra of $\tM$ through the
identification $\lambda^{\omega\cdot E}(t)=\lambda^{\omega}(t)$,
where $E$ is the minimal expectation from $M$ onto $N$ and $\omega$ is a
faithful normal semifinite weight on $N$.
We denote by $\tE$ the extension of $E$ to $\tM$ leaving
$\lambda^{\varphi_0\cdot E}(t)$ fixed, which is a normal conditional
expectation from $\tM$ onto $\tN$ \cite{KL}.
Even when $M\supset N$ has an infinite index, if $E$ is a unique normal
conditional expectation from $M$ onto $N$, we use the same convention.

Before introducing the canonical extension, we need some preparation.

\begin{definition} Let $M$ be an infinite factor.
A pair $(\varphi, \rho)$ consisting of a faithful normal semifinite weight
$\varphi$ on $M$ and $\rho\in \eM$ is said to be an \textit{invariant pair}
if the following hold:
$$\varphi\cdot\rho=d(\rho)\varphi,$$
$$\varphi\cdot E_\rho=\varphi.$$
\end{definition}

Note that the above conditions are equivalent to
$d(\rho)\varphi\cdot \phi_\rho=\varphi,$
thanks to $\rho\cdot \phi_\rho=E_\rho$ and
$\phi_\rho\cdot\rho=\id$.

\begin{lemma} Let $M$ be a type $\mathrm{III}$ factor.
Then, the following hold:
\begin{itemize}
\item[$(1)$] If $\varphi$ is a faithful normal semifinite weight on $M$
and $\rho\in \eM$, then $\varphi\cdot \phi_{\rho}$ is a faithful normal
semifinite weight on $M$.
\item[$(2)$] If $(\varphi,\rho)$ is an invariant pair, then
$\rho$ commutes with $\sigma^\varphi_t$ for $t\in \R$.
\item[$(3)$] Let $\psi$ be a dominant weight
(see \cite{CT} for the definition).
Then, for every sector in $\sM$, there exists a representative
$\rho\in \eM$ such that $(\psi, \rho)$ is an invariant pair.
\item[$(4)$] Let $\psi$ be a dominant weight on $M$ and $\rho\in \eM$.
If $(\psi,\rho)$ is an invariant pair, then there exist irreducible
$\rho_i\in \eM$ and isometries $V_i\in (\rho_i, \rho)$ $i=1,2,\cdots, n$
such that each $V_i$ belongs to the centralizer $M_\psi$ and
$$\rho(x)=\sum_{i=1}^nV_i\rho_i(x)V_i^*,\quad x\in M.$$
Moreover, each $(\psi,\rho_i)$ is an invariant pair.
\item[$(5)$] Let $\psi$ be a dominant weight and $\rho_1,\rho_2\in \eM$.
If $(\psi,\rho_1)$ and $(\psi,\rho_2)$ are invariant pairs,
then $(\rho_1,\rho_2)$ is in the centralizer $M_\psi$.
\end{itemize}
\end{lemma}

\begin{proof} (1). Since $\varphi\cdot \rho^{-1}$ is a faithful normal
semifinite
weight on $\rho(M)$, so is
$\varphi\cdot \phi_{\rho}=\varphi\cdot \rho^{-1}\cdot E_\rho$ on $M$.

(2). Let $N=\rho(M)$ and $\varphi_0$ be the restriction of
$\varphi$ to $N$, which is a faithful normal semifinite weight on $N$
by assumption.
Then, we have
$$\varphi_0\cdot \rho=d(\rho)\varphi,$$
$$\varphi_0\cdot E_\rho=\varphi.$$
Regarding $\rho$ as an isomorphism from $M$ to $N$, we get
$$\rho\cdot \sigma^\varphi_t=\sigma^{\varphi_0}_t\cdot\rho
=\sigma^\varphi_t\cdot\rho.$$

(3). This is \cite[Lemma 2.12]{ILP}.

(4). Let $N=\rho(M)$ and $\psi_0$ be the restriction of $\psi$ to $N$.
Then, we have common continuous decomposition \cite{Loi}, \cite{KL}
$$M=M_\psi\rtimes_{\theta_0}\R \supset N
=N_{\psi_0}\rtimes_{\theta_0}\R,$$
such that
$$M\cap N'=(M_\psi\cap N_{\psi_0}')^{\theta_0}.$$
Therefore, each non-zero projection $p\in M\cap N'$ is equivalent to
1 in $M_\psi$ (more strongly in $M_\psi^{\theta_0}$), and there exist
isometries $V_i\in M_\psi$, $i=1,2,\cdots, n$
such that each $p_i=V_iV_i^*$ is a minimal projection in $M\cap N'$
and $\{p_i\}_{i=1}^n$ forms a partition of unity.
Setting $\rho_i(x)=V_i^*\rho(x)V_i$, we get the first assertion.
Since the minimal expectation for $p_iMp_i\supset Np_i$
is given by $E_\rho(x)p_i/E_\rho(p_i),$ we have
$$E_{\rho_i}(x)
=\frac{d(\rho)}{d(\rho_i)}V_i^* \rho \cdot \phi_\rho(V_ixV_i^*)V_i.
=\frac{d(\rho)}{d(\rho_i)}\rho_i\cdot\phi_\rho(V_ixV_i^*),
\quad x\in M,$$
thanks to the local index formula for the minimal conditional
expectation \cite{H} and the fact that $\mathrm{Ad}(V_i^*)$ is
an isomorphism from $p_iMp_i\supset N_{p_i}$ to $M\supset \rho_i(M)$.
This implies
$$\phi_{\rho_i}(x)=\frac{d(\rho)}{d(\rho_i)}\phi_\rho(V_ixV_i^*),
\quad x\in M.$$
Thus,
$$d(\rho_i)\psi\cdot \phi_{\rho_i}=d(\rho)\psi\cdot
\phi_\rho(V_i\cdot V_i^*)=\varphi(V_i\cdot V_i^*).$$
Since $V_i$ is in the centralizer of $\psi$, this implies
that $(\psi,\rho_i)$ is an invariant pair.

(5). Thanks to (4), the general case is reduced to the case where
$\rho_1$ and $\rho_2$ are irreducible, and we make this assumption.
If $(\rho_1,\rho_2)=0$, we have nothing to prove, and so we assume
$\rho_1=\mathrm{Ad}(u)\cdot\rho_2$ for some unitary $u\in M$.
Since $\sigma^\psi_t$ commutes with $\rho_1$ and $\rho_2$,
$\sigma^\psi_t(u)$ is proportional to $u$, and there exists a scalar
$c>0$ such that $\sigma^\psi_t(u)=c^{it}u$.
However, the KMS condition and $\psi\cdot \rho_1=\psi\cdot\rho_2$ imply
$c=1$, which shows $u\in M_\psi$.
\end{proof}

Now we introduce the canonical extension of an endomorphism.

\begin{theorem} Let $M$ be a type $\mathrm{III}$ factor and $\varphi$
be a faithful normal semifinite weight on $M$.
Then, for every $\rho\in \eM$, there exists a unique endomorphism
$\tr\in \mathrm{End}(\tM)$ satisfying
$$\tr(x)=\rho(x), \quad x\in M$$
$$\tr(\lambda^\varphi(t))=d(\rho)^{it}
[D\varphi\cdot \phi_\rho:D\varphi]_t\lambda^\varphi(t).$$
Under the identification $\lambda^\psi(t)
=[D\psi:D\varphi]_t\lambda^\varphi(t)$, $\tr$ does not
depend on the choice of the weight $\varphi$.
\end{theorem}

\begin{proof} Let $\psi$ be a dominant weight.
First we assume that $(\psi,\rho)$ is an invariant pair, and
$N$ and $\psi_0$ are as in the proof of Lemma 2.3, (4).
Thanks to Lemma 2.3, (2), we may consider $\rho$ an isomorphism
from $M$ onto $N$ intertwining $\sigma^{\psi}$ and $\sigma^{\psi_0}$.
Therefore, $\rho$ extends to an isomorphism $\tr$ from $\tM$ onto
$\tN$ sending $\lambda^\psi(t)$ to $\lambda^{\psi_0}(t)$.
Since $\psi$ satisfies $\psi=\psi_0\cdot E_\rho$, we may regard $\tN$
as a subalgebra of $\tM$ identifying $\lambda^{\psi_0}(t)$ with
$\lambda^\psi(t)$.
We show that $\tr$ satisfies the required property.
Indeed,
$$\tr(\lambda^\varphi(t))=\tr([D\varphi:D\psi]_t\lambda^\psi(t))
=\rho([D\varphi:D\psi]_t)\lambda^\psi(t).$$
Regarding $\rho$ as an isomorphism from $M$ onto $N$ again, we have
$$\rho([D\varphi:D\psi]_t)=[D\varphi\cdot\rho^{-1}:D\psi\cdot\rho^{-1}]_t
=[D\varphi\cdot\rho^{-1}\cdot E_\rho:D\psi\cdot\rho^{-1}\cdot E_\rho]_t,$$
and so
$$\tr(\lambda^\varphi(t))=[D\varphi\cdot\phi_\rho:D\psi\cdot\phi_\rho]_t
\lambda^\psi(t)
=d(\rho)^{it}[D\varphi\cdot \phi_\rho:D\varphi]_t\lambda^\varphi(t).$$

Now, we treat the general case.
Thanks to Lemma 2.3, (3), it suffices to show the statement for
$\rho_1=\mathrm{Ad}(u)\cdot \rho$, where $u$ is a unitary in $M$,
and $\rho$ and $\psi$ are as before.
Indeed,
\begin{eqnarray*}\mathrm{Ad}(u)\cdot \tr(\lambda^\varphi(t))
&=&d(\rho)^{it}u\sigma^{\varphi\cdot \phi_\rho}_t(u^*)
\lambda^{\varphi\cdot \phi_\rho}(t)\\
&=&d(\rho)^{it}[D\varphi\cdot \phi_\rho\cdot \mathrm{Ad}(u^*):
D\varphi\cdot \phi_\rho]_t\lambda^{\varphi\cdot \phi_\rho}(t)\\
&=&d(\rho)^{it}[D\varphi\cdot\phi_{\rho_1}:D\varphi]_t\lambda^\varphi(t).
\end{eqnarray*}
Therefore, $\tr_1:=\mathrm{Ad}(u)\cdot \tr$ has the desired property.
\end{proof}

The same formula of the Connes cocycle as above also appears in
\cite{L5}.

We call the above $\tr$ the \textit{canonical extension} of $\rho$.
When, $\alpha$ is an automorphism, the left inverse $\phi_\alpha$ is
nothing but the inverse of $\alpha$.
Therefore, our definition of the canonical extension generalizes that
in the automorphism case.

Let $M\supset N$ be properly infinite von Neumann algebras.
Then, the standard representation of $M$ is automatically standard for $N$
and the product of two modular conjugation $J_NJ_M$ makes sense.
The canonical endomorphism $\gamma\in  \mathrm{End}(M)$ for the inclusion
$M\supset N$ is defined by $\gamma(x)=J_NJ_M x J_MJ_N\in N$ \cite{L1}.
In general, modular conjugations depend on the choice of natural cones.
However, a different choice of modular conjugations amounts to only a
perturbation of $\gamma$ by an inner automorphism of $N$.
Therefore, we call any endomorphism in $\mathrm{End}(M)$ of the form
$\mathrm{Ad}(u)\cdot \gamma$ with a unitary $u\in N$ the canonical
endomorphism as well.

\begin{proposition} Let $M$ be a type $\mathrm{III}$ factor and
$\rho,\mu\in \eM$.
Then, the following hold:
\begin{itemize}
\item [$(1)$] $(\rho,\mu)\subset (\tr,\tm)$.
\item [$(2)$] $\widetilde{\rho\cdot\mu}=\tr\cdot\tm$.
\item [$(3)$] Let $N$ be a subfactor of $M$ with a finite index and $\gamma$
be the canonical endomorphism for $M\supset N$.
Then, $\widetilde{\gamma}$ is the canonical endomorphism for $\tM\supset
\tN$.
\item [$(4)$] $\tau\cdot \tr=d(\rho)\tau$.
\item [$(5)$] $\tr\cdot \theta_t=\theta_t\cdot \tr$.
\end{itemize}
\end{proposition}

\begin{proof} In the proof of Theorem 2.4, we have already shown that
for $\rho\in \eM$ and a unitary $u\in M$, we have
$$\widetilde{(\mathrm{Ad}(u)\cdot \rho)}=\mathrm{Ad}(u)\cdot \tr.$$
Thus, in order to prove the statements we may replace $\rho$ and $\mu$ with
unitary equivalent endomorphisms.
Thanks to Lemma 2.3, (3), we may and do assume that $(\psi,\rho)$ and
$(\psi,\mu)$ are invariant pairs, where $\psi$ is a dominant weight.
Therefore, we have
$$\tr(\lambda^\psi(t))=\tm(\lambda^\psi(t))=\lambda^\psi(t).$$
In this situation, (2), (4), and (5) are obvious and (1) follows from
Lemma 2.3, (5).

To prove (3), we consider the common continuous decomposition of
$M\supset N$.
We may and do assume that $\psi$ is of the form $\psi_0\cdot E$,
where $\psi_0$ is a dominant weight on $N$ and $E$ is the minimal
expectation from $M$ to $N$.
Then, we have the common continuous decomposition
$$M=M_\psi\rtimes_{\theta_0}\R \supset
N=N_{\psi_0}\rtimes_{\theta_0}\R .$$
Let $\tau_0$ be the trace whose dual weight is $\psi$.
Then, thanks to \cite[Section 4]{L2}, the canonical endomorphism of
$M_\psi\supset N_{\psi_0}$ scales $\tau_0$ by $[M:N]_0$.
Let $J_{M_\psi}$ and $J_{N_{\psi_0}}$  be the modular conjugations of
$M_\psi$ and $N_{\psi_0}$, and $\gamma_0$ be the corresponding canonical
endomorphism of $M_\psi\supset N_{\psi_0}$.
To define the canonical endomorphism of $M\supset N$ and $\tM\supset \tN$,
we utilize the modular objects of these algebras naturally coming from
the modular objects of $M_\psi$ and $N_{\psi_0}$ through the crossed
products
as computed in \cite[Lemma 2.8]{Ha1}.
Let $\gamma$ and $\gamma_1$ be such canonical endomorphisms of
$M\supset N$ and $\tM\supset \tN$ respectively.
Then, $(\psi,\gamma)$ is an invariant pair \cite[Section 3]{I1}.
On the other hand, applying \cite[Lemma 3.1]{I1} twice and using
\cite[Lemma 2.8]{Ha1} for the computation of
$\Delta_\psi$ and $\Delta_{\psi_0}$, we know that $\gamma_1$ is an
extension of $\gamma$ leaving $\lambda^\psi(t)$ invariant.
Therefore, $\gamma_1$ coincides with the canonical extension of $\gamma$.
\end{proof}

\section{Modular endomorphisms}
Let $M$ be a type III factor and $\alpha$ be an automorphism of $M$.
In \cite[Proposition 5.4]{HS2}, it was shown that $\ta$ is inner
if and only if $\alpha$ is the composition of an inner automorphism
with an extended modular automorphism.
We adopt this property as the definition of an endomorphism counterpart
of an extended modular automorphism.

\begin{definition} Let $M$ be a type III factor and $\rho\in \eM$.
We say that $\rho$ is a \textit{modular endomorphism} if $\tr$ is
an inner endomorphism, that is, there exist isometries
$\{V_i\}_{i=1}^n\subset \tM$ with mutually orthogonal ranges and
$\sum_{i=1}^nV_iV_i^*=1$ such that
$$\tr(x)=\sum_{i=1}^nV_ixV_i^*,\quad x\in \tM. $$
The system of isometry $\{V_i\}_{i=1}^n$ satisfying the above condition
is called an \textit{implementing system} for $\tr$.
\end{definition}

We denote by $\eMm$ the set of modular endomorphisms of $M$.
Thanks to Proposition 2.5, (1), an endomorphism unitary equivalent to
a modular endomorphism is again a modular endomorphism,
and so it makes sense to call a sector modular if its representatives
are modular endomorphisms.
We denote by $\sMm$ the set of sectors of modular endomorphisms.
It is also easy to see from Proposition 2.5, (1), (2) that $\eMm$
is closed under forming a direct sum and a product of finitely many
elements.
What is not clear for the moment is whether $\eMm$ is closed
under conjugation and irreducible decomposition,
which will be shown later.
The number $n$ of the isometries of the implementing system in
Definition 3.1 is nothing but $d(\rho)$, which also will be shown
later.

Let $M$ be as above and $Z(\tM)$ be the center of $\tM$.
Then, the restriction of $\theta$ to $Z(\tM)$ is an ergodic action of the
real number group $\R$.
Thus, by Mackey's point realization theorem \cite{M1}, there exist
a standard Borel space $X_M$, a probability measure $\mu_M$ on $X_M$,
and an ergodic flow $F^M$ on $(X_M,\mu_M)$ such that
$Z(\tM)=L^\infty(X_M,\mu_M)$
and $\theta_t(f)(\omega)=f(F^M_{-t}\omega)$ for $f\in L^\infty(X_M,\mu_M)$,
$\omega\in X_M$, $t\in \R$.
(For simplicity, we often omit $F^M$ and denotes $F^M_{-t}\omega$ just by
$\omega\cdot t$).
$F^M$ is called  the smooth flow of weights \cite{CT}.

Let $K$ be a compact group.
A Borel map $c:X_M\times \R\longrightarrow K$ is said to be a
cocycle of
$F^M$ if for fixed $s,t\in \R$, the cocycle relation
$$c(\omega,s)c(\omega\cdot s,t)=c(\omega,s+t)$$
holds for almost all $\omega\in X_M$.
Two cocycles agree on the outside of a null set are identified as usual.
We denote by $Z^1(F^M,K)$ the set of all $K$-valued cocycles of $F^M$.
Two cocycles $c$ and $c'$ are said to be equivalent
(or cohomologous) if there exists a Borel map $a:X_M\longrightarrow K$
such that for fixed $t\in \R$
$$c(\omega,t)=a(\omega)c'(\omega,t)a(\omega\cdot t)^{-1}$$
holds for almost all $\omega\in X_M$.
A cocycle equivalent to a constant function $e\in K$ is said to be a
coboundary.
We denote by $H^1(F^M,K)$ the set of equivalence classes of $K$-valued
cocycles of $F^M$.
For the significance of cocycles from the viewpoint of Mackey's notion
of virtual groups, readers are refered to \cite{M3} and \cite{Z1}.

As Connes and Takesaki showed that the group of extended modular
automorphisms divided by the inner automorphism group is isomorphic to
$H^1(F^M, U(1))$ \cite{CT}, we show that $\sMm$ is ``isomorphic" to
the union $\bigcup_{n=1}^\infty H^1(F^M, U(n))$ of the unitary group
$U(n)$-valued cohomology classes.

Let $\rho\in \eMm$ and $\{V_i\}_{i=1}^n$ be an implementing system for
$\tr$.
Since $(\id,\tr)$ is globally preserved by $\theta$,
$c(t)_{ij}:= V_i^*\theta_t(V_j)$ belongs to $Z(\tM)$.
$c(t)=(c(t)_{ij})$ is regarded as a $U(n)$-valued Borel function on
$X_M\times \R$.

\begin{lemma} Let $\rho$ and $c(t)$ be as above.
Then,
\begin{itemize}
\item[$(1)$] $c$ is a $U(n)$-valued cocycle of $F^M$, whose cohomology
class depends only on the sector of $\rho$.
\item[$(2)$] Let $c'$ be a cocycle equivalent to $c$.
Then, there exists an implementing system $\{V'_i\}_{i=1}^n$ for
$\tr$ such that $c'(t)_{ij}={V'_i}^*\theta_t(V'_j)$.
\end{itemize}
\end{lemma}

\begin{proof} (1). The cocycle relation of $c$ follows from
$$\sum_{j=1}^nc(s)_{ij}\theta_s(c(t)_{jk})=\sum_{j=1}^nV_i^*
\theta_s(V_jV_j^*\theta_{t}(V_k))=c(s+t)_{ik}.$$
Let $\{W_i\}_{i=1}^m\subset \tM$ be another implementing system for $\tr$,
and set $a_{ij}=V_i^*W_j$.
Since $a(t)=(a(t)_{ij})$ is a matrix valued function that is unitary,
we have $n=m$.
Setting  $b(t)_{ij}=W_i^*\theta_t(W_j)$, we get
$$c(t)_{ij}=V_i^*\theta_t(V_j)=\sum_{k=1,l}^nV_i^*W_kW_k^*
\theta_t(W_lW_l^*V_j)
=\sum_{k=1,l}^na_{ik}b(t)_{kl}\theta_t(a_{jl}^*),$$
which shows that $b$ and $c$ are equivalent.
Let $u$ be a unitary of $M$ and $\rho_1=\mathrm{Ad}(u)\cdot \rho$.
Then, $\{uV_i\}_{i=1}^n$ is an implementing system for $\tr_1$.
Since $u$ is fixed by $\theta_t$, $\rho_1$ gives the same cohomology class
as $\rho$.

(2). We take a $U(n)$-valued function $f$ satisfying
$$c(t)_{ij}=\sum_{k,l=1}^n f_{ik}c'(t)_{kl}\theta_t(f_{jl}^*).$$
Then, $V'_i:=\sum_{j=1}^n f_{ji}V_j$ has a desired property.
\end{proof}

We introduce a map
$$\delta_m: \sMm \longrightarrow \bigcup_{n=1}^\infty H^1(F^M, U(n))$$
sending $[\rho]$ to $[c]$ in the above lemma.

As described in \cite[Section 2]{Z1}, $Z^1(F^M,U(n))$ has natural
three operations in analogous to the unitary representation theory of
compact groups: direct sum $c\oplus c'$, tensor product $c\otimes c'$ and
complex conjugate $\bc$.
$H^1(F^M,U(n))$ inherits these operations and we use the same notation for
the cohomology classes as well.

\begin{theorem} Let $\delta_m$ be as above and $\rho_1,\rho_2\in \eMm$.
Then,
\begin{itemize}
\item[$(1)$] $\delta_m$ is a bijection.
\item[$(2)$]  We have
$\delta_m([\rho_1]\oplus [\rho_2])=\delta_m([\rho_1])\oplus
\delta_m([\rho_2]).$
\item[$(3)$] We have
$\delta_m([\rho_1][\rho_2])=\delta_m([\rho_1])\otimes
\delta_m([\rho_2]).$
\item[$(4)$] $[\rho_1]$ is the conjugate sector of $[\rho_2]$ if and only if
$\delta_m([\rho_1])=\overline{\delta_m([\rho_2])}$.
In particular, $\sMm$ is closed under conjugation.
\item[$(5)$] $\delta_m$ is grade preserving in the sense that
$\delta_m([\rho])\in H^1(F^M, U(n))$ if and only if  $d(\rho)=n$.
\end{itemize}
\end{theorem}

\begin{proof} (1). First we show that $\delta_m$ is injective.
Let $\rho_1,\rho_2\in \eMm$ satisfying
$\delta_m([\rho_1])=\delta_m([\rho_2])$.
Then, thanks to Lemma 3.2, (2), there exist implementing systems
$\{V^{(1)}_i\}_{i=1}^n$ and $\{V^{(2)}_i\}_{i=1}^n$ for
$\tr_1$ and $\tr_2$ respectively such that they give the same cocycle.
We set $u=\sum_{i=1}^n V^{(1)}_i{V^{(2)}_i}^*$, which is a unitary in
$\tM^\theta=M$ satisfying $\rho_1=\mathrm{Ad}(u)\cdot \rho_2$.
This shows that $\delta_m$ is injective.

Next we show that $\delta_m$ is surjective.
Let $\psi$ be a dominant weight on $M$ and
$M=M_\psi\rtimes_{\theta_0}\R$ be the continuous decomposition.
We take the trace $\tau_0$ on $M_\psi$ whose dual weight is $\psi$,
and take the implementing one-parameter unitary group
$\{u(s)\}_{s\in \R}$ for $\theta_0$.
We assume that $M_\psi$ acts on a Hilbert space $H$.
Then, the Takesaki duality theorem \cite[Chapter I]{NT} implies that there
exists an
isomorphism $\Phi$ from $\tM$ to $M_\psi\otimes B(L^2(\R))$ such
that
$$\Phi(x)=\pi_{\theta_0}(x), \quad x\in M_\psi,$$
$$\Phi(u(s))=1\otimes l(s),$$
$$\Phi(\lambda^\psi(t))=1\otimes m(t),$$
where $l(s)$ is the left regular representation of
$\R$ and $m(t)$ is the multiplication operator of $e^{-its}$.
Moreover, $\theta_t$ is identified with
$(\theta_{0,t}\otimes \mathrm{Ad}(r(t)))$ under $\Phi$, or more precisely
we have  $\Phi\cdot \theta_t=(\theta_{0,t}\otimes \mathrm{Ad}(r(t)))
\cdot \Phi$.
Using $\Phi$, we identify the flow $(Z(\tM),\theta)$ with
$$(Z(M_\psi\otimes B(L^2(\R))), \theta_0\otimes \mathrm{Ad}(r))
\cong (Z(M_\psi),\theta_0).$$

Let $c\in Z^1(F^M,U(n))$ be a given cocycle, where $L^\infty(X_M,\mu_M)$
is understood as the point realization of $Z(M_\psi)$.
We take a system of isometries $\{W_i\}_{i=1}^n$ in $M_\psi$ satisfying
$\sum_{i=1}^nW_iW_i^*=1$ and set
$$W(t)=\sum_{i,j=1}^n c(t)_{ij}W_i\theta_{0,t}(W_j^*)\in M_\psi.$$
Then, $W(t)$ is a unitary satisfying
\begin{eqnarray*}W(s+t)&=&
\sum_{i,j,k=1}^n c(s)_{ik}\theta_{0,s}(c(t)_{kj})W_i\theta_{0,s+t}(W_j^*)\\
&=&\sum_{i,j,k=1}^n c(s)_{ik}W_i
\theta_{0,s}(W_k^*c(t)_{lj}W_l\theta_{0,t}(W_j^*))\\
&=&W(s)\theta_{0,s}(W(t)),
\end{eqnarray*}
which means that $W(t)$ is a $\theta_0$-cocycle.
Since $\theta_0$ is stable \cite[Theorem III.5.1]{CT},
every $\theta_0$-cocycle is
a coboundary and there exists a unitary $w$ such that
$W(t)=w^*\theta_{0,t}(w)$.
Let $V_i=wW_i$.
Then, $V_i$ is an isometry in $M_\psi$ satisfying $\sum_{i=1}^nV_iV_i^*=1$
and
$$\theta_{0,t}(V_i)=\sum_{j=1}^n c(t)_{ji}V_j.$$
For $x\in M_\psi$ we set,
$$\rho_0(x)=\sum_{i=1}^nV_ixV_i^*,$$
$$\phi_0(x)=\frac{1}{n}\sum_{i=1}^nV_i^*xV_i,$$
and $E_0=\rho_0\cdot \phi_0$.
Then, $\rho_0$ is an endomorphism of $M_\psi$ commuting with
$\theta_{0,t}$, $\phi_0$ is a left inverse of $\rho_0$ satisfying
$n\tau_0\cdot\phi_0=\tau_0$, and $E_0$ is a $\tau_0$-preserving
conditional expectation from $M_\psi$ onto the image of $\rho_0$.
We introduce $\rho\in \mathrm{End}(M)$ that is an extension of
$\rho_0$ to $M$ leaving $u(t)$ invariant.
Indeed, such $\rho$ exists because the following hold:
$$\Phi(\rho_0(x))=\sum_{i=1}^n(V_i\otimes 1)\Phi(x)(V_i^*\otimes 1),
\quad x\in M_\psi,$$
$$\Phi(u(t))=\sum_{i=1}^n(V_i\otimes 1)\Phi(u(t))(V_i^*\otimes 1).$$
In a similar way, we define a left inverse $\phi$ of $\rho$ and a
conditional expectation $E$ from $M$ onto the image of $\rho$ by
the relations
$$\Phi(\phi(x))=\frac{1}{n}\sum_{i=1}^n(V_i^*\otimes 1)\Phi(x)
(V_i\otimes 1), \quad x\in M,$$
and $E=\rho\cdot \phi$.
Note that $\rho$, $\phi$, and $E$ are extensions of $\rho_0$, $\phi_0$,
and $E_0$ respectively.

We claim that $[M:\rho(M)]_0=n^2$, $\phi=\phi_\rho$, $E=E_\rho$, and
$(\psi,\rho)$ is an invariant pair.
Indeed, it suffices to show $\mathrm{Index}\; E=[M:\rho(M)]_0=n^2$,
and the other claims can be shown easily.
We define $T\in M_\psi$ by
$$T={\sqrt{n}}\sum_{i=1}^n V_iV_i.$$
Then, for $x\in M_\psi$ and $u(t)$, we have
\begin{eqnarray*}
T^*E(Txu(t))&=&T^*E_0(Tx)u(t)
=\frac{1}{n}\sum_{i,j=1}^nT^*V_iV_j^*TxV_jV_i^*u(t)\\
&=&\sum_{i,j=1}^nV_i^*V_jxV_jV_i^*u(t)
=xu(t),
\end{eqnarray*}
which shows $\{T\}$ is a one-element Pimsner-Popa basis and
$\mathrm{Index}\; E=T^*T=n^2$ \cite{PP}, \cite{Wat}.
To show that $E$ is minimal, it suffices to show
$T^*xT=n^2E(x)$ for all $x\in M\cap \rho(M)'$
\cite[Theorem 2.12.3]{Wat}.
Thanks to the Connes-Takesaki relative commutant theorem
\cite[Theorem II.5.1]{CT},
for $x\in M\cap \rho(M_\psi)'$ we have
$$x_{ij}:=V_i^*xV_j\in M\cap M_\psi'=Z(M_\psi),$$
and
$$x=\sum_{i,j=1}^nV_ix_{ij}V_j^*.$$
Thus,
$$T^*xT=n\sum_{i=1}^nx_{ii}=n^2E(x),$$
which shows the claims.

Since $(\psi,\rho)$ is an invariant pair, $\tr$ is an extension of
$\rho$ leaving $\lambda^\psi(t)$ invariant.
Therefore, we have
$$\Phi(\tr(x))=\sum_{i=1}^n(V_i\otimes 1)\Phi(x)(V_i^*\otimes 1),
\quad x\in \tM,$$
which shows that $\rho$ is a modular endomorphism carrying $c$.

(2) and (3) are easy, and (5) follows from the above proof.

(4). Assume that $\delta_m([\rho_1])=[c]$ and
$\delta_m([\rho_2])=[\bc]$.
Let $\{V_i\}_{i=1}^n$ and $\{W_i\}_{i=1}^n$ be implementing
systems for $\tr_1$ and $\tr_2$ satisfying
$$\theta_t(V_i)=\sum_{j=1}^n c(t)_{ji}V_j,$$
$$\theta_t(W_i)=\sum_{j=1}^n \bc(t)_{ji}W_j.$$
We define two isometries $R$ and $S$ by
$$R=\frac{1}{\sqrt{n}}\sum_{i=1}^n V_iW_i,$$
$$S=\frac{1}{\sqrt{n}}\sum_{i=1}^n W_iV_i.$$
$R$ and $S$ belong to $\tM^\theta=M$ satisfying
$$R\in (\id,\rho_1\cdot\rho_2)$$
$$S\in (\id,\rho_2\cdot\rho_1)$$
$$S^*\rho_2(R)=R^*\rho_1(S)=\frac{1}{n}.$$
Thus, $[\rho_1]=\overline{[\rho_2]}$ \cite{L3}.
Since $\delta_m$ is a bijection, this also shows that
$[\rho]\in \sMm$ implies $\overline{[\rho]}\in \sMm$ and
$\delta_m(\overline{[\rho]})=\overline{\delta_m([\rho])}$.
\end{proof}

The following is a useful criterion for an irreducible
endomorphism to be modular.

\begin{proposition} Let $M$ be a type $\mathrm{III}$ factor.
\begin{itemize}
\item[$(1)$] Let $\rho\in \eM$ be an irreducible endomorphism.
If there exists a non-zero element in $(\id,\tr)$, $\rho$ is a
modular endomorphism.
\item[$(2)$]  Every irreducible component of a modular endomorphism
is again a modular endomorphism.
In particular $\sMm$ is closed under irreducible decomposition.
\end{itemize}
\end{proposition}

\begin{proof} (1). First, we show that $(\id,\tr)$ contains at least one
isometry if $(\id,\tr)\neq \{0\}$.
Let $\mathcal{PI}(\id,\tr)$ be the set of partial isometries in
$(\id,\tr)$.
Using the polar decomposition, we know that $\mathcal{PI}(\id,\tr)$ is not
empty.
We introduce an order into $\mathcal{PI}(\id,\tr)$ as follows:
For two element $v_1, v_2\in \mathcal{PI}(\id,\tr)$, we say that $v_2$
dominates $v_1$ if the following holds
$$v_1=v_2v_1^*v_1.$$
Note that when this is the case, $v_1^*v_1\in Z(\tM)$, and
$v_2^*v_1=v_1^*v_1$, $v_2v_1^*=v_1v_1^*$ hold.
We show that $\mathcal{PI}(\id,\tr)$ is an inductively ordered set
with this order.
Let $\mathcal{T}$ be a totally ordered subset of $\mathcal{PI}(\id,\tr)$.
Since $\{vv^*\}_{v\in \mathcal{T}}$ and $\{v^*v\}_{v\in \mathcal{T}}$ are
increasing nets of projections in $Z(\tM)$ and $\tM$, they converge
to projections in $Z(\tM)$ and in $\tM$ respectively in strong topology.
We claim that the net $\{v\}_{v\in \mathcal{T}}$ converges in strong $*$
topology.
Let $v_1,v_2 \in \mathcal{T}$ with $v_2$ dominating $v_1$, and $\xi$ be
a vector in the Hilbert space that $\tM$ acts on.
Then,
$$||(v_2-v_1)\xi||^2=\inpr{(v_2^*v_2-v_1^*v_1)\xi}{\xi},$$
$$||(v_2^*-v_1^*)\xi||^2=\inpr{(v_2v_2^*-v_1v_1^*)\xi}{\xi}.$$
Therefore, the strong $*$ limit exists, and it gives a majorant of
$\mathcal{T}$, which shows that $\mathcal{PI}(\id,\tr)$ is
inductively ordered.
Now, we apply Zorn's lemma to $\mathcal{PI}(\id,\tr)$ and take a maximal
element $V\in \mathcal{PI}(\id,\tr)$.
We show that $V$ is an isometry.
Suppose $V^*V\neq 1$.
Since $V^*V\in Z(\tM)$ and $\theta$ acts on $Z(\tM)$ ergodically,
there exists some $t\in \R$ such that
$\theta_t(V^*V)(1-V^*V)\neq 0$.
Therefore, $V_1:=V+\theta_t(V)(1-V^*V)\in \mathcal{PI}(\id,\tr)$
dominates $V$, which is contradiction.
Thus, $V$ is an isometry.

Let $\{V_\lambda\}_{\lambda\in \Lambda}$ be a maximal set of isometries
in $(\id, \tr)$ with mutually orthogonal ranges
(such a set exists thanks to Zorn's lemma again.), and
let $P=\sum_{\lambda\in \Lambda}V_\lambda V_\lambda^*$.
If $P=1$, we are done, and so we assume $P\neq 1$.
Let $\mathcal{PI}_{P^\bot}(\id,\tr)$ be the subset of elements in
$\mathcal{PI}(\id,\tr)$ with range projections orthogonal to $P$.
We claim that there exists a non-zero element in
$\mathcal{PI}_{P^\bot}(\id,\tr)$.
Indeed, since $(\tr,\tr)^\theta=(\rho,\rho)=\C 1$, there exists
some $t\in \R$ such that $(1-P)\theta_t(P)\neq 0$, which shows that
there exists non-zero $x\in (\id,\tr)$ with $(1-P)x\neq 0$.
Using the polar decomposition, we get the claim.
We introduce an order into $\mathcal{PI}_{P^\bot}(\id,\tr)$ as before
and take a maximal element $W$.
Since $\{V_\lambda\}_{\lambda\in \Lambda}$ is a maximal set, $W$ is not
an isometry, and we set $W^*W=z\neq 1$.
We claim $(1-P)(1-z)=0$.
Indeed, by maximality of $W$, $(1-P)(1-z)x=(1-P)x(1-z)=0$  holds for
$x\in (\id,\tr)$, and $\theta_t((1-P)(1-z))x=0$ for $x\in (\id,\tr)$
as well for all $t\in \R$ because $(\id,\tr)$ is globally
invariant under $\theta$.
Thus, we get
$$(\bigvee_{t\in \R}\theta_t((1-P)(1-z)))x=0, \quad x\in
(\id,\tr).$$
Since $(\tr,\tr)^\theta=\C 1$, this implies $(1-P)(1-z)=0$.

Let $V'_\lambda=V_\lambda(1-z)$.
Then, $\{V'_\lambda\}_{\lambda\in \Lambda}$ have mutually orthogonal ranges
such that
$$\sum_{\lambda\in \Lambda}V'_\lambda {V'_\lambda}^*=P(1-z)=(1-z).$$
Thus, we get
$$(\tM (1-z)\supset \tr(\tM)(1-z))\cong ((B(\ell^2(\Lambda))\otimes
\tM(1-z))\supset
(\C 1\otimes \tM(1-z))).$$
Since there exists a conditional expectation from $\tM$ to $\tr(\tM)$
satisfying the Pimsner-Popa inequality \cite{HK3}, \cite{PP},
$\Lambda$ is a finite set, and we identify $\Lambda$ with $\{1,2,\cdots,
n\}$.
We take $t\in \R$ satisfying $z_1:=\theta_t(z)(1-z)\neq 0$, and
set $\Omega\subset X_M$ to be the Borel subset corresponding $z_1$.
Let $R_i=V'_i z_1$, $S_i=\theta_t(V_i)z_1$, $i=1,2,\cdots n$, and
$S_{n+1}=\theta_t(W)z_1$.
Then,  $\{R_i\}_{i=1}^{n}\subset (\id,\tr)$ and
$\{S_i\}_{i=1}^{n+1}\subset (\id,\tr)$ satisfy
$$R_i^*R_j=\delta_{i,j}z_1,\quad i,j=1,2,\cdots n,$$
$$\sum_{i=1}^nR_iR_i^*=z_1,$$
$$S_i^*S_j=\delta_{i,j}z_1,\quad i,j=1,2,\cdots,n+1,$$
$$\sum_{i=1}^{n+1}S_iS_i^*\leq z_1.$$
We set $f=(f_{ij})=(R_i^*S_j)$, which is regarded as an $n$ by $n+1$
matrix-valued
function on $\Omega$.
However $f$ satisfies
$$\sum_{i=1}^n f_{ij}^*f_{ik}=\delta_{j,k}z_1,\quad j,k=1,2,
\cdots, n+1,$$
which is contradiction.
Thus, we conclude $P=1$, and
$$\tr(x)=\sum_{\lambda\in \Lambda}V_\lambda xV_\lambda^*.$$

(2). Let $\rho$ be a modular endomorphism with an implementing system
$\{V_i\}_{i=1}^n$ for $\tr$, and $\sigma$ be an irreducible component of
$\rho$.
We take a non-zero element $T\in (\sigma,\rho)$.
Thanks to Proposition 2.5, (1),  $T^*V_i$ belongs to
$(\id, \widetilde{\sigma})$, $i=1,2,\cdots, n$.
Since $T^*=\sum_{i=1}^n T^* V_iV_i^*$, there exists some $i$ such that
$T^*V_i$ is non-zero.
Thus, $\sigma$ is a modular endomorphism thanks to (1).
\end{proof}

Though it is shown in the above that $\sMm$ is closed under irreducible
decomposition, we have not given a criterion for a modular endomorphism
to be irreducible, or that of how to decompose it when reducible, in
terms of ergodic theory yet.
We discuss this issue now.

Let $X$ be a standard Borel space, $\mu$ be a probability measure on $X$,
and $G$ be a locally compact group ergodically acting on $(X,\mu)$
as a non-singular transformation group.
For a compact group $K$, we define a $K$-valued cocycle $c$ as
in the case of the flow.
We denote by $K_c$ the closed subgroup generated by the image of $c$.
We collect necessary results on cocycles from R. Zimmer's fundamental
paper \cite[Part I. Section 3]{Z1} in the next theorem in order to introduce
a few notions for cocycles.
Though Zimmer treated only the measure preserving case in \cite{Z1},
his argument works for the non-singular case with a formal modification
(see also \cite{F}):

\begin{theorem}[R. Zimmer]
Let the notations be as above, and $c\in Z^1(G,K)$.
Then, the following hold:
\begin{itemize}
\item[$(1)$] There exists a closed subgroup $H\subset K$ such that
$c$ is equivalent to an $H$-valued cocycle $c'$, and $c'$ is never
equivalent to $c''$ with $K_{c''}$ a proper subgroup of $H$.
$H$ is uniquely determined up to conjugacy.
We call such $c'$ a minimal cocycle and call $H$ the minimal subgroup
of $c$.
\item[$(2)$] Let $K_1$ be a compact group containing $K$.
Then, if $c$ is a minimal cocycle as a $K$-valued cocycle,
it is the case as a $K_1$-valued cocycle.
\item[$(3)$] Let $c$ be a minimal $K$-valued cocycle such that
$K=K_c$, and $\pi$ be a continuous homomorphism from $K$ to a compact
group $K_1$.
Then, $\pi\cdot c$ is a minimal cocycle.
\item[$(4)$] Let $c$ be a minimal $K$-valued cocycle with $K=K_c$,
and $(\pi_1,H_{\pi_1})$, $(\pi_2,H_{\pi_2})$
be finite dimensional unitary representations of $K$.
If $a$ is a Borel map from $X$ to $\mathrm{Hom}(H_1,H_2)$
such that for every fixed $k$,
$$a(\omega)\pi_2(c(\omega,k))=\pi_1(c(\omega,k))a(\omega\cdot k),$$
holds for almost all $\omega\in X$,
then $a$ is constant almost everywhere with the value in $(\pi_1,\pi_2)$.
\end{itemize}
\end{theorem}

We say that a closed subgroup $K\subset U(n)$ is irreducible if
the defining representation of $K$ on $\C^n$ is irreducible.

\begin{corollary} Let $M$ be a type $\mathrm{III}$ factor and $\rho$ be
a modular endomorphism with $\delta_m([\rho])=[c]$.
Then, $\rho$ is irreducible if and only if the minimal subgroup of $c$ is
irreducible.
\end{corollary}

\begin{proof} Thanks to theorem 3.3, it suffices to show that the minimal
subgroup of $c$ is irreducible if and only if $c$ is never equivalent to
a direct sum of two cocycles.
Assume that the minimal subgroup of $c$ is irreducible.
If $c$ were equivalent to $c_1\oplus c_2$, $c_i\in Z^1(G,U(n_i))$,
$i=1,2$, with $n_1+n_2=n$,
$c$ could be considered as a $U(n_1)\times U(n_2)$-valued cocycle,
and the minimal subgroup of $c$ would be conjugate to a subgroup
of $U(n_1)\times U(n_2)$ in $U(n)$ thanks to Theorem 3.5, (2).
This is contradiction and $c$ is never equivalent to a direct sum of
two cocycles.
The converse also follows from Theorem 3.5, (2).
\end{proof}

We end this section with giving an answer to the problem
about the type II and type III principal graphs described in
Introduction.

Let $M\supset N$ be an inclusion of type III factors of a finite index with
the minimal expectation $E$.
We say that $M$ and $N$ has the common flow of weights if
$Z(\tN)=Z(\tM)$.
When this is the case, we have common central decomposition
$$(\tM\supset \tN)=
\int_{X_M}^{\oplus}(\tM(\omega)\supset \tN(\omega))d \mu_M(\omega),$$
and we could define the type II principal graph for that of
$\tM(\omega)\supset \tN(\omega)$ for $\omega$ in a conull set.
However, it would be cumbersome to treat ``measurable field of
subfactors'' (though it should not be too hard to do so,
c.f. \cite{Wi}.)
Instead, we formulate the problem using a global term.
Let
$$N\subset M\subset M_1\subset \cdots M_n\subset \cdots$$
be the tower for $N\subset M$.
Then, using the minimal conditional expectation in each step,
we have the ``tower for the core inclusion'' \cite{KL}:
$$\tN\subset \tM\subset \tM_1\subset \cdots \tM_n\subset \cdots.$$
We still denote by $\theta$ the natural extension of $\theta$
to $\tM_n$ leaving the Jones projections fixed.
Then, we have
$$M_n\cap M'=(\tM_n\cap \tM')^\theta. $$
We say that \textit{graph change occurs} if
$Z(\tM)\vee (M_n\cap M')$ does not coincides with
$\tM_n\cap \tM'$ for some $n$.
Let $E_n$ be the minimal expectation from $M_n$ onto
$M_{n-1}$ and $\tE_n$ be the natural extension of $E_n$ to
$\tM_{n}$, which is a conditional expectation from $\tM_{n}$
onto $\tM_{n-1}$.
Then, using $E_{n+1}$ and $\tE_{n+1}$, we can show that
$Z(\tM)\vee (M_n\cap M')\neq \tM_n\cap \tM'$ implies
$Z(\tM)\vee (M_{n+1}\cap M')\neq \tM_{n+1}\cap \tM'$.
In a similar way using downward basic construction and mirroring
$J_M\cdot J_M$, $J_{\tM}\cdot J_{\tM}$, we can also show that
graph change occurs for $M\supset N$ if and only if
it occurs for $M_1\supset M$.

A modular endomorphism is said to be non-trivial if it is
irreducible and not equivalent to identity.

\begin{theorem} Let $M\supset N$ be an inclusion of type
$\mathrm{III}_0$ factors with a finite index and with the common
flow of weights.
Then, graph change occurs for $M\supset N$ if and only if
a non-trivial modular endomorphism appears in some power of the
canonical endomorphism $\gamma$ for $M\supset N$.
\end{theorem}

\begin{proof} Assume that graph change occurs for $N\subset M$.
Taking sufficiently large $n$, we may assume
$Z(\tM)\vee (M_{2n}\cap M')\neq \tM_{2n}\cap \tM'$.
Let $\gamma$ be the canonical endomorphism for $M\supset N$.
Then, thanks to Proposition 2.5, this is equivalent to
$$(\widetilde{\gamma}^n,\widetilde{\gamma}^n)\neq
(\gamma^n,\gamma^n)\vee Z(\tM).$$
Therefore, either of the following two holds:
(i) $\gamma^n$ contains two irreducibles $\rho_1,\rho_2\in \eM$
such that $[\rho_1]\neq [\rho_2]$ and $(\tr_1,\tr_2)\neq \{0\}$.
(ii) $\gamma^n$ contains an irreducible $\rho\in \eM$ such that
$(\tr,\tr)\neq Z(\tM)$.
Indeed, let $\{p_i\}_{i=1}^n$ be a partition of unity consisting
of minimal projections in $M\cap \gamma^n(M)'$.
If $$X\in (\widetilde{\gamma}^n,\widetilde{\gamma}^n)\setminus
(\gamma^n,\gamma^n)\vee Z(\tM),$$
there exist some $i$ and $j$ such that
$$p_iXp_j \in (\widetilde{\gamma}^n,\widetilde{\gamma}^n)\setminus
(\gamma^n,\gamma^n)\vee Z(\tM).$$
If the irreducible components corresponding to $p_i$ and $p_j$ are
different, (i) occurs.
If they are the same, (ii) occurs.

In the case (i), the Frobenius reciprocity implies
$(\id,\widetilde{\rho_2\cdot\br_1})\neq \{0\}$ (though $\tM$ is not
a factor, the argument in \cite[Section 2]{I3} works).
Thanks to Lemma 3.3, this shows that $\rho_2\cdot\br_1$,
which is contained in $\gamma^{2n}$, contains a modular endomorphism.
Since $\rho_2\cdot\br_1$ does not contain identity, it is indeed
a non-trivial one.

In the case (ii), we take $X\in (\tr,\tr)\setminus Z(\tM)$.
Let $R\in (\id,\br\cdot \rho)$, $\bR\in (\id,\rho\cdot\br)$
be isometries satisfying the usual relation \cite{L3}
$$R^*\br(\bR)=\bR^*\rho(R)=\frac{1}{d(\rho)}.$$
Then, $X\bR\in (\id,\widetilde{\rho\cdot\br})$.
Let $\sigma$ be an irreducible sector contained in $\rho\cdot\br$
that is not equivalent to identity, and $\{Y(\sigma)_i\}_{i=1}^k$
be a basis of $(\sigma,\rho\cdot\br)$.
If $Y(\sigma)_i^*X\bR\in (\id,\widetilde{\sigma})$ is not zero for some
$i$, we are done thanks to Lemma 3.3.
Thus, we assume $Y(\sigma)_i^*X\bR=0$ for every $\sigma$ not equivalent
to identity.
This would imply that $\bR\bR^*X\bR=X\bR$.
However, since $\rho(R)$ commutes with $X$, we would get
$$X=d(\rho)\rho(R^*)X\bR=d(\rho)\rho(R^*)\bR\bR^*X\bR=\bR^*X\bR\in
Z(\tM),$$
which is contradiction.
Thus, $\rho\cdot\br$ contains a non-trivial modular endomorphism,
and so does $\gamma^{2n}$.

The converse can be shown easily.
\end{proof}

\begin{remark}
For a type $\mathrm{III}_\lambda$, $0<\lambda\leq 1$ factor, the flow is
a transitive and every $U(n)$-valued cocycle is equivalent to a direct sum
of
$U(1)$-valued  cocycles that come from homomorphisms from
the stabilizer subgroup of a point \cite{M3}, $\Z$ for $\lambda\neq
1$
and $\R$ for $\lambda=1$.
This means that every modular endomorphism is decomposed into usual
modular automorphisms.
Therefore, the above theorem generalizes the previous one obtained in
\cite[Theorem 3.5]{I1}.
\end{remark}

There exist plenty of examples in the type $\mathrm{III}_0$ case,
where higher dimensional non-trivial modular endomorphisms
appear in the canonical endomorphisms (see Section 5 and Proposition A.5 ).

\section{Connes-Takesaki modules}
For a type III factor $M$, we denote by $\mathrm{Aut}(F^M)$ the
set of automorphisms of $Z(\tM)$ that commute with the restriction of
$\theta$ to $Z(\tM)$.
For an automorphism $\alpha$ of $M$, the restriction of $\ta$ to
$Z(\tM)$ belongs to $\mathrm{Aut}(F^M)$, which is called
the Connes-Takesaki module of $\alpha$.
This correspondence gives a homomorphism
$$\mathrm{\mo}: \mathrm{Aut}(M) \longrightarrow \mathrm{Aut}(F^M),$$
which is called the fundamental
homomorphism.
(Though the original definition of Connes and Takesaki looks different
from ours \cite{CT}, this is an equivalent description of it.)
The purpose of this section is to introduce the Connes-Takesaki
module for endomorphisms.

\begin{definition} Let $M$ be a type III factor and $\rho\in \eM$.
We say that \textit{$\rho$ has a Connes-Takesaki module} if
$M\supset \rho(M)$ has the common flow of weights.
When $\rho$ has a Connes-Takesaki module, we denote by
$\mo(\rho)$ the restriction of $\tr$ to $Z(\tM)$.
$\mo(\rho)\in \mathrm{Aut}(F^M)$ is called the
\textit{Connes-Takesaki module of $\rho$}.
\end{definition}

We denote by $\eMC$ (respectively $\eMt$) the set of endomorphisms in
$\eM$ with Connes-Takesaki modules (respectively trivial
Connes-Takesaki modules).
Note that $\mo(\rho)$ depends only on the sector of $\rho$.
We denote by $\sMC$  and $\sMt$ the corresponding subsets of $\sM$.

\begin{proposition} Let $M$ be a type $\mathrm{III}$ factor. Then,
\begin{itemize}
\item[$(1)$] If $\rho_1,\rho_2\in \eM$ have Connes-Takesaki modules,
then so does $\rho_1\cdot \rho_2$ and
$\mo(\rho_1\cdot\rho_2)=\mo(\rho_1)\cdot\mo(\rho_2)$.
\item[$(2)$] Let $\rho, \rho_i\in \eM$, $i=1,2,\cdots, n,$ satisfying
$[\rho]=\oplus_{i=1}^n[\rho_i]$.
Then, $\rho$ has a Connes-Takesaki module if and only if
each $\rho_i$ does and $\mo(\rho_i)$ does not depend on $i$.
When it is the case, $\mo(\rho)=\mo(\rho_i)$.
\item[$(3)$] $\rho\in \eM$ has a Connes-Takesaki module if and only
if the conjugate $\br$ does.
When it is the case, $\mo(\br)=\mo(\rho)^{-1}$.
\item[$(4)$] Every modular endomorphism has a trivial Connes-Takesaki
module.
\item[$(5)$] Let $N\subset M$ be a subfactor with a finite index.
If $M$ and $N$ have the common flow of weights, then the canonical
endomorphism $\gamma$ for $M\supset N$ has a trivial Connes-Takesaki
module.
\end{itemize}
\end{proposition}
\begin{proof} (1). This is trivial.

(2). We take isometries $V_i\in (\rho_i,\rho)$ satisfying
$$\rho(x)=\sum_{i=1}^n V_i\rho_i(x)V_i^*,\quad x\in M.$$
Assume $\rho\in \eMC$.
Thanks to Lemma 2.4, (1), we have
$$\tr_i(z)=V_i^*\tr(z)V_i=\tr(z),\quad z\in Z(\tM),$$
which shows that $\rho_i\in \eMC$ for every $i$ and
$\mo(\rho_i)$ does not depend on $i$.
Next we assume that $\rho_i\in \eMC$ for every $i$ and
$\mo(\rho_i)$ does not depend on $i$, which is denoted by $\nu$.
Then,
$$\tr(z)=\sum_{i=1}^nV_i\tr_i(z)V_i^*=\sum_{i=1}^n \nu(z)V_iV_i^*=\nu(z),
\quad z\in Z(\tM),$$
which shows $\rho\in \eMC$.

(3). Since $M\supset \br(M)$ is the dual inclusion of $M\supset \rho(M)$,
the first part is obvious.
The second part follows from (1) and (2) because $\rho\cdot\br$ contains
identity.

(4). This follows from the definition of modular endomorphisms.

(5). This follows from Proposition 2.5, (3).
\end{proof}

Let $M\supset N$ be an inclusion of type $\mathrm{III}_\lambda$,
$0<\lambda<1$ factors with a finite index, and $E$ be the minimal
expectation from $M$ onto $N$.
We say that $M\supset N$ has common discrete decomposition if
there exist an inclusion of type $\mathrm{II}_\infty$ factors
$P\supset Q$ with a finite index, a minimal expectation $E_0$ from
$P$ to $Q$ preserving the trace, and a trace scaling automorphism
$\theta_0$ of $P$ globally preserving $Q$ such that
$$(M\supset N)\cong (P\rtimes_{\theta_0}\Z\supset
Q\rtimes_{\theta_0}\Z).$$
Note that in \cite{Loi}, common discrete decomposition with respect
to any faithful normal expectation is discussed, while here
we consider the minimal expectation only.

The following lemma is standard and we omit the proof.

\begin{lemma} Let $M\supset N$ be an inclusion of type
$\mathrm{III}_\lambda$, $0<\lambda<1$ factors with a finite index,
and $E$ be the minimal expectation from $M$ to $N$.
Then, the following are equivalent:
\begin{itemize}
\item[(1)] $M\supset N$ has common discrete decomposition.
\item[(2)] $M\supset N$ has the common flow of weights.
\item[(3)] For some (and any) generalized trace $\psi_0$ on $N$
(see \cite{C1} for the definition, it is also called $\lambda$-trace
in \cite{St}), $\psi_0\cdot E$ is a generalized trace on $M$.
\end{itemize}
\end{lemma}

\begin{corollary} Let $M$ be a type $\mathrm{III}_\lambda$, $0<\lambda<1$,
factor, and $\rho\in \eM$.
Then,
\begin{itemize}
\item[(1)] If $\rho$ is irreducible, $\rho$ has a Connes-Takesaki module.
\item[(2)] $M\supset \rho(M)$ has common discrete decomposition
if and only if $\rho$ has a Connes-Takesaki module,  which is further
equivalent to that every irreducible component of $\rho$ has the
same Connes-Takesaki module.
\end{itemize}
\end{corollary}

\begin{proof} (1). This follows from the condition (3) of Lemma 4.3.
(2). This follows from Proposition 4.2, (2) and Lemma 4.3.
\end{proof}

Let $M$ be a type $\mathrm{III}_{\lambda}$ factor, $0<\lambda<1$,
and $\psi$ be a generalized trace on $M$.
We set $T=-2\pi/\log \lambda$.
Then, the modular automorphism group $\{\sigma^\psi_t\}_{t\in \R}$
has period $T$ and we consider $\sigma^\psi$
an action of $\T:=\R/T\R$.
We denote by $\hM$ the crossed product
$M\rtimes_{\sigma^\psi}\T$, which is a type $\mathrm{II}_\infty$
factor generated by a copy of $M$ and the implementing
unitary representation $\{\lambda_d^\psi(t)\}$ of $\T$.
If $\varphi$ is another weight on $M$ satisfying
$[D\varphi:D\psi]_T=1$, $\{[D\varphi:D\psi]_t\}$ is a
$\sigma^\psi$-cocycle as a $\T$ action, and
we identify $[D\varphi:D\psi]_t\lambda_d^\psi(t)$ with
$\lambda_d^\varphi(t)$.
We denote by $\tau_d$ and $\theta_d$ the natural trace on $\hM$ and
the dual action of $\sigma^\psi$ respectively.
In analogous to the canonical extension, we can construct a functor
from $\eMt$ to $\mathrm{End}(\hM)_0$.
Note that thanks to Proposition 4.2, $\eMt$ is a full subcategory of $\eM$.

\begin{proposition} Let $M$ be a type $\mathrm{III}_{\lambda}$ factor
and $\psi$ be a generalized trace.
Then, for every $\rho\in \eMt$, there exists a unique endomorphism
$\hr\in \mathrm{End}(\hM)$ satisfying
$$\hr(x)=\rho(x), \quad x\in M$$
$$\hr(\lambda_d^\psi(t))=d(\rho)^{it}
[D\psi\cdot \phi_\rho:D\psi]_t\lambda_d^\psi(t).$$
Moreover, for $\rho,\mu\in \eMt$, the following hold:
\begin{itemize}
\item [$(1)$] $(\rho,\mu)\subset (\hr,\widehat{\mu})$.
\item [$(2)$] $\widehat{\rho\cdot\mu}=\hr\cdot\widehat{\mu}$.
\item [$(3)$] Let $N$ be a subfactor of $M$ with a finite index and $\gamma$
be the canonical endomorphism for $M\supset N$.
If $M\supset N$ has common discrete decomposition,
then, $\widehat{\gamma}$ is the canonical endomorphism for
$\hM\supset \widehat{N}$.
\item [$(4)$] $\tau_d\cdot \hr=d(\rho)\tau_d$.
\item [$(5)$] $\hr\cdot \theta_d=\theta_d\cdot \hr$.
\end{itemize}
\end{proposition}

\begin{proof} Since the center of $\tM$ is generated by $\lambda^\psi(T)$,
we have $\tr(\lambda^\psi(T))=\lambda^\psi(T)$ by assumption,
which implies
$$[Dd(\rho)\psi\cdot\phi_\rho:D\psi]_T=1. $$
Therefore, there exists a unitary $u\in M$ such that
$d(\rho)\psi\cdot\phi_\rho=\psi\cdot \mathrm{Ad}(u)$
thanks to \cite[Th{\'e}or{\`e}me 4.3.2]{C1}.
This means that $(\psi,\mathrm{Ad}(u)\cdot \rho)$ is an invariant pair.
The rest of the proof is the same as those of Theorem 2.4 and
Proposition 2.5 if the dominant weight is replaced with the generalized
trace.
\end{proof}

We end this section with a remark on the type $\mathrm{II}_\infty$ case.

\begin{remark}
Let $M$ be a type $\mathrm{II}_\infty$ factor and $\tau$ be a faithful
normal
semifinite trace on $M$.
Then, the canonical extension and the Connes-Takesaki modules
make sense for $M$ as well.
We introduce a scalar valued module $\Mo(\rho)$ of
$\rho\in \eM$ by
$$\tau\cdot \rho=d(\rho)\Mo(\rho)\tau.$$
Note that every inclusion of type $\mathrm{II}_\infty$ factors
with a finite index comes from the tensor product of a common type I factor
and an inclusion of $\mathrm{II}_1$ factors.
Therefore, the restriction of $\tau$ to the image of $\rho$ is a
semifinite trace and the above definition makes sense.
The center of $\tM$ is generated by $\lambda^\tau(t)$, and we have
$$\tr(\lambda^\tau(t))=d(\rho)^{it}
(\frac{d\tau\cdot \phi_\rho}{d\tau})^{it}\lambda^\tau(t)
=\Mo(\rho)^{-it}(\frac{d\tau\cdot E_\rho}{d\tau})^{it}\lambda^\tau(t).$$
This means that $\rho$ has a Connes-Takesaki module if and only if
$E_\rho$ is trace preserving.
When this is the case, we have
$$\mo(\rho)(\lambda^\tau(t))=\Mo(\rho)^{-it}\lambda^\tau(t).$$
Therefore, Proposition 4.2, (2) implies that $M\supset \rho(M)$ is an
extremal inclusion (see \cite{PP} for the definition) if and only if
each irreducible component of $\rho$ has the same $\Mo$.
\end{remark}

\section{Minimal Actions of Compact Groups}
In this section, we investigate the structure of minimal actions of
compact groups on type III factors applying our machinery to the
corresponding Roberts actions.
Modular endomorphisms and the Connes-Takesaki modules for
endomorphisms provide new invariants for minimal actions of
compact groups through Roberts actions.

\subsection{Minimal actions as dual actions}

Let $K$ be a compact group and $\mathrm{Rep}(K)$ be the tensor
category of finite dimensional unitary representations of $K$.
Roughly speaking, a Roberts action of the dual of $K$ on a von Neumann
algebra $N$ is a functor from a (sufficiently big) subcategory $\cR$ of
$\mathrm{Rep}(K)$, (called a ring in \cite{NT}), to $\eN$,
where the set of arrows from $\rho_1\in  \eN$ to $\rho_2\in \eN$ is
$(\rho_1,\rho_2)$.
For the precise definition of the Roberts action and
the crossed product by it, we refer to \cite{R} and \cite{NT}
(see also \cite{DR1}).
For each equivalence class of irreducible representations of $K$,
we choose a representative $(\pi, H_\pi)$, and denote by $\hK$
the collection of them.
We always assume that $\cR$ contains every member of $\hK$ and its
complex conjugate representation.
We often omit to specify $\cR$ in a Roberts action,
and say ``a Roberts action of $\hK$" if there is no possibility of
confusion (although it is a bit sloppy terminology).
For $(\pi_1,H_{\pi_1}), (\pi_2,H_{\pi_2})\in \cR$, we denote by
$(\pi_1,\pi_2)$ the set of $K$-homomorphisms from $H_{\pi_1}$ to
$H_{\pi_2}$, that is, the set of arrows from $\pi_1$ to $\pi_2$.

Let $M$ be a factor and $\alpha$ be an action of $K$ on $M$.
$\alpha$ is said to be \textit{minimal} if $\alpha$ is faithful
and the relative commutant $M\cap {M^\alpha}'$ of the fixed point
subalgebra $M^\alpha$ is trivial \cite{PW}.
Note that when $\alpha$ is minimal in our sense, the crossed product
$M\rtimes_\alpha K$ is automatically a factor isomorphic to
$M^\alpha\otimes B(L^2(K))$.
Indeed, to show this, we may assume that $M^\alpha$ is infinite by taking
tensor product with a type I factor if necessary.
Then, a similar argument as in \cite[Lemma III.3.4]{AHKT} implies that
for each irreducible representation $\pi$ of $K$, there exists a Hilbert
space
$\cH_\pi$ in $M$ that is globally invariant under $\alpha$ such that the
restriction of
$\alpha$ to $\cH_\pi$ is equivalent to $\pi$.
Therefore, \cite[Proposition 6.9]{R} implies the claim.

More strongly, \cite[Theorem 6.5]{R} shows that $M$ is the crossed product
$N\rtimes_\beta \widehat{K}$ by some Roberts action $\beta$,
and that $\alpha$ is its dual action $\widehat{\beta}$, where
$N=M^\alpha$.
The relationship between $\beta_\pi$ and $\cH_\pi$ is as follows:
Let $\{W(\pi)\}_{i=1}^{d(\pi)}$ be an orthonormal basis for $\cH_\pi$,
where $d(\pi)$ is the dimension of $\pi$.
Then, $\beta_\pi\in \eN$ is given by
$$\beta_\pi(x)=\sum_{i=1}^{d(\pi)}W(\pi)_ixW(\pi)_i^*,\quad x\in N.$$
On the other hand, $\cH_\pi$ is recovered from $\beta_\pi$ via
$$\cH_\pi=\{W\in M;\; Wx=\beta_\pi(x)W,\quad x\in N \}.$$
Note that $\beta$ is uniquely determined up to equivalence defined in
\cite[Definition 5.6]{R}.

We will use these notations throughout this section whenever a minimal
action $\alpha$ of a compact group $K$ on a factor $M$ is given.

\begin{lemma} Let $N$ be an infinite factor, and $\beta$ be a Roberts action
of $\hK$ on $N$.
Then, the dual action of $\beta$ is minimal if and only if
$\beta_\pi$ is irreducible for every irreducible representation
$(\pi,H_\pi)\in \hK$.
\end{lemma}

\begin{proof} We set $M=N\rtimes_\beta \hK$, $\alpha=\widehat{\beta}$,
and
$$E(x)=\int_K\alpha_k(x)dm_K(k), \quad x\in M,$$
where $m_K$ is the Haar measure of $K$.
Let $\{W(\pi)_i\}_{i=1}^{d(\pi)}$ be as above.
For $x\in M$, we define $x(\pi)_i\in N$ by
$$x(\pi)_i=d(\pi)E(W(\pi)_ix).$$
Then, $x$ has the following formal expansion (see \cite[Section 3]{ILP}):
$$x``="\sum_{(\pi,H_\pi)\in \hK}\sum_{i=1}^{d(\pi)}W(\pi)_i^*x(\pi)_i.$$
More precisely, the above summation converges in the GNS Hilbert space
topology with respect to an $\alpha$-invariant normal state, and
$\{x(\pi)_i\}$ completely determines $x$.
Using this expansion, we can see that $x\in M\cap N'$ if and only if
$x(\pi)_i\in (\id,\beta_\pi)$.
Therefore, $\alpha$ is minimal if and only if $\beta_\pi$ does not contain
identity for every non-trivial irreducible representation $\pi$.
However, the latter is equivalent to the statement that $\beta_\pi$ is
irreducible for very irreducible $\pi$ thanks to the Frobenius reciprocity
for the usual compact group representations and that of sectors \cite{I4}.
\end{proof}

Let $\{u(k)\}_{k\in K}$ be an $\alpha$-cocycle.
We denote by $\alpha^u$ the perturbation
$\mathrm{Ad}(u(k))\cdot \alpha_k$  of $\alpha$ by the cocycle $u$.
$\alpha$ is called stable if every $\alpha$-cocycle is a coboundary.

The following statements are probably well-known among specialists,
and actually some of them already exist in the literature.
However, since we cannot find complete proofs in the literature,
we provide them here for convenience of readers.

\begin{proposition} Let $\alpha$ be a minimal action of a compact
group $K$ on a factor $M$.
We denote by $E$ the unique normal conditional expectation from
$M$ onto $M^\alpha$ obtained by taking average over $K$.
Then, the following hold:
\begin{itemize}
\item[$(1)$] Every cocycle perturbation of $\alpha$ is again minimal.
\item[$(2)$] If $M$ is of type $\mathrm{II}_1$, so is $M^\alpha$.
$\alpha$ is stable in this case \cite[Theorem 12]{Was1}.
\item[$(3)$] If $M$ is of type $\mathrm{II}_\infty$, so is $M^\alpha$.
$\alpha$ is stable in this case.
\item[$(4)$] If $M$ is of type $\mathrm{III}_0$, $M^\alpha$ is of type
$\mathrm{III}$.
$\alpha$ is stable in this case
(more generally, $\alpha$ is stable whenever $M^\alpha$ is of type
$\mathrm{III}$).
\item[$(5)$] $M$ is of type $\mathrm{III}_\lambda$, $\lambda\neq 0$
and $M^\alpha$ is of type $\mathrm{II}$, if and only if
the center of $\alpha_K$ contains the modular automorphism group
$\{\sigma^{\varphi\cdot E}_t\}_{t\in \R}$, for some
faithful normal semifinite weight $\varphi$ on $M^\alpha$
(which is actually a trace if this is the case).
$\alpha$ is not stable in this case.
\end{itemize}
\end{proposition}

\begin{proof}(1). First we claim that the second dual action
$\alpha\otimes \Ad(r)$ on $M\otimes B(L^2(K))$ is minimal whenever
$\alpha$ is so.
To show the claim, we may assume that  $M^\alpha$ is infinite as usual.
Thus, we have a Hilbert space $\cH_\pi$ in $M$ for each
$(\pi, H_\pi)\in \hK$ as before.
This implies that there exists a Hilbert space $\cH$ in $M$ such that
the restriction of $\alpha$ to $\cH$ is equivalent to the regular
representation.
Therefore, $(M,\alpha)$ is conjugate to
$(M\otimes B(L^2(K)), \alpha\otimes \Ad(r))$, which shows the claim.
Let $\alpha^u$ be a cocycle perturbation of $\alpha$.
Since $M\rtimes_\alpha K\subset M\otimes B(L^2(K))$ is
an irreducible inclusion and there exists a projection
$p\in M\rtimes_\alpha K$ such that
$$(M^{\alpha^u}\subset M)\cong
(p(M\rtimes_\alpha K)p\subset p(M\otimes B(L^2(K)))p),$$
which shows that $\alpha^u$ is minimal.

(2),(3). The first part is easy.
Stability in these two cases (and the case (4) as well) follows from
(1) and Connes' 2 by 2 matrix trick \cite{C5}.

(4). We show that if $M^\alpha$ is of type II,
$M$ is either of type II or of type $\mathrm{III}_\lambda$,
$\lambda\neq 0$.
Let $\tau$ be a faithful normal semifinite trace on $M^\alpha$.
Then, the centralizer $M_{\tau\cdot E}$ is an intermediate subfactor
between $M^\alpha$ and $M$, and in particular
$M_{\tau\cdot E}$ is a factor.
Therefore, the Connes spectrum $\Gamma(\sigma^{\tau\cdot E})$
coincides with the Arveson spectrum
$\mathrm{Sp}(\sigma^{\tau\cdot E})$  \cite{C1}, \cite{St}.
Thus, $M$ is either of type II or type $\mathrm{III}_\lambda$,
$\lambda\neq 0$ depending on the period of the modular
automorphism group $\{\sigma^{\tau\cdot E}_t\}_{t\in \R}$.

(5). Note that part of the proof has been already done in the above.
Assume that $\alpha_K$ contains
$\{\sigma^{\varphi\cdot E}_t\}_{t\in \R}$, for some
faithful normal semifinite weight $\varphi$ on $M^\alpha$.
Then, the centralizer $M_{\varphi\cdot E} $ is an intermediate
subfactor of $M^\alpha\subset M$.
Let $F$ be the restriction of $E$ to $M_{\varphi\cdot E}$,
which is a normal conditional expectation from $M_{\varphi\cdot E}$
to $M^\alpha$.
Since the restriction of $\varphi\cdot E$ to $M_{\varphi\cdot E}$
is nothing but $\varphi\cdot F$ and it is semifinite,
$M_{\varphi\cdot E}$ is of type II, and so is $M^\alpha$ as well because
$\sigma^{\varphi\cdot F}$ is trivial.
This finishes the proof of the first statement.

Now, we show that $\alpha$ is not stable in the above situation with
type III $M$ .
First we assume that $M^\alpha$ is of type $\mathrm{II}_\infty$.
Let $p\in M^\alpha$ be a non-zero finite projection, and
$V\in M$ be an isometry satisfying $VV^*=p$.
We set $u(k)=V^*\alpha_k(V)$, which is an $\alpha$-cocycle.
It is easy to show that $M^{\alpha^u}$ is finite and $u$ is not a
coboundary.
When $M^\alpha$ is of type $\mathrm{II}_1$, a similar construction
identifying $M$ with a corner of $M\otimes B(H)$ works.
\end{proof}

\begin{remark}
In the above situation, assume that $M^\alpha$ is of type
$\mathrm{II}_\infty$
with a faithful normal trace $\tau$.
Then, the type of $M$
(or the period of $\{\sigma^{\tau\cdot E}\}_{t\in \R}$)
is completely determined by $\Mo(\beta_\pi)$ for
the Roberts action $\beta$ satisfying $M=M^\alpha\rtimes_\beta \hK$,
where $\Mo$ is the scalar Connes-Takesaki module introduced
in Remark 4.6.
Indeed, as in \cite[Section 3]{ILP}, we can obtain the action of
the modular automorphism group on $\cH_\pi$ as follows:
$$\sigma^{\tau\cdot E}_t(W)=\Mo(\beta_\pi)^{it}W,
\quad W\in \cH_\pi.$$
\end{remark}

It is possible to show that the canonical extension gives continuous
homomorphism from $\mathrm{Aut}(M)$ to $\mathrm{Aut}(\tM)$ in the
$u$-topologies in the same spirit of the Cones-Takesaki's proof of
continuity of the fundamental homomorphism \cite{CT}.
However, we just mention here that when $\alpha$ is an action of a locally
compact group $G$ on $M$ with a faithful normal invariant state, then it is
particularly easy to show that $G\ni g\mapsto \ta_g$ is continuous.
On the other hand, when $\beta$ is a Roberts action of $\hK$ on a factor
$N$,
Proposition 2.5 shows that $\beta$  extends to $\tb$ on $\tN$ via the
canonical extension.

\begin{lemma} Let $\alpha$ be an action of a compact group $K$ on a factor
$M$,
and $E$ be the normal conditional expectation from $M$ onto $N:=M^\alpha$
given
by the average of $\alpha_k$ over $K$.
We regard $\tN$ as a subalgebra of $\tM$ identifying
$\lambda^{\varphi_0}(t)$
with $\lambda^{\varphi_0\cdot E}(t)$, where $\varphi_0$ is a faithful normal
semifinite weight on $\tN$.
Then,
\begin{itemize}
\item [$(1)$] The fixed point subalgebra of $\tM$ under $\ta$ coincides with
$\tN$.
\item [$(2)$] When $\alpha$ is a minimal action that is the dual action of a
Roberts action $\beta$ of $\hK$ on $N$.
Then, $\tM$ is naturally isomorphic to $\tN\rtimes_{\tb}\hK$.
\end{itemize}
\end{lemma}

\begin{proof} (1). We set $\varphi=\varphi_0\cdot E$.
Let $\tE$ be the normal conditional expectation from $\tM$ onto the fixed
point
subalgebra of $\tM$ under $\ta$ obtained by average of $\ta_k$ over $K$.
Since $\ta_k$ leaves $\lambda^\varphi(t)$ fixed for all
$k\in K$ and $t\in \R$, the image of $\tE$ coincides with $\tN$.

(2). Let $\psi_0$ be a dominant weight on $N$, and $\psi=\psi_0\cdot E$.
For each $\beta_\pi$, $(\pi,H_\pi)\in \hK$, we may and do assume that
$(\psi_0,\beta_\pi)$ is an invariant pair thanks to Lemma 2.3, (3).
Then, the modular automorphism group $\{\sigma^\psi_t\}_{t\in \R}$
acts on $\cH_\pi$ trivially \cite[Section 3]{ILP} and
$\tb_\pi$ leaves $\lambda^{\psi_0}(t)$ invariant.
Therefore, $\tb_\pi$ is given by
$$\tb_\pi(x)=\sum_{i=1}^{d(\pi)}W(\pi)_i x W(\pi)_i^*,\quad x\in \tN.$$
Thus, (1) and \cite[Theorem 6.5]{R} imply that $\tM=\tN\rtimes_{\tb}\hK$
and $\ta$ is the dual action of $\tb$.
\end{proof}

\subsection{A Roberts action of modular endomorphisms}
The purpose of this subsection is to show how a system of modular
endomorphisms gives rise to a Roberts action, and to see the structure of
its dual
action.
There are two options to do so, (of course they are essentially
the same), and we start with the one using the minimal subgroup
of a cocycle.

Let $N$ be a type III factor and $\cD\subset \eNm$ be a
countable family of modular endomorphisms of $N$.
By adding new endomorphisms to $\cD$ if necessary, we may and do assume
the following conditions:
\begin{itemize}
\item[(i)] $\cD$ is closed under taking a product of any two members
in $\cD$.
\item[(ii)] For every pair $\rho_1,\rho_2\in \cD$ there exists
$\rho\in \cD$ such that $[\rho_1]\oplus [\rho_2]=[\rho]$.
\item[(iii)] If $\sigma\in \eN$ is contained in $\rho\in \cD$,
there exists an endomorphism in $\cD$ that is equivalent to $\sigma$.
\item[(iv)] Every $\rho\in \cD$ has conjugate (up to equivalence)
in $\cD$.
\end{itemize}

For each $\rho\in \cD$ we choose $c_\rho\in Z^1(F^N,U(d(\rho)))$
satisfying $\delta_m([\rho])=[c_\rho]$, and set
$$K_1=\prod_{\rho\in \cD}U(d(\rho)),$$
$$c=(c_\rho)_{\rho\in \cD}\in Z^1(F^N,K_1),$$
where $U(n)$ is the unitary group.
Thanks to Theorem 3.5, we may and do assume that $c$ is already a minimal
cocycle.
We denote by $K$ the minimal subgroup of $c$.

Let $(\pi_\rho, H_{\pi_\rho})$ be the representation of $K$
on $H_{\pi_\rho}=\C^{d(\rho)}$ that is the projection
onto the $\rho$-component.
Then, we have $c_\rho=\pi_\rho\cdot c$.
We denote by $\cR$ the collection of $(\pi_\rho, H_{\pi_\rho})$,
$\rho\in \cD$.

Now, we construct a Roberts action $\beta$ of $\cR$ on $N$
sending $\pi_\rho$ to $\rho$.
Thanks to Lemma 3.2, (2), we can choose an implementing system
$\{V(\rho)_i\}_{i=1}^{d(\rho)}$ for $\tr$ giving $c_\rho$.
Let $\cK_\rho$ be the linear span of $\{V(\rho)_i\}_{i=1}^{d(\rho)}$,
which is a Hilbert space in $\tN$.
For $\rho_1,\rho_2\in \cD$, we denote by $\cK_{\rho_2}\cK_{\rho_1}^*$
the linear span of elements of the form $V_2V_1^*$, $V_1\in \cK_1$,
$V_2\in \cK_2$, which is naturally identified with
$\mathrm{Hom}(\cK_{\rho_1},\cK_{\rho_2})$ and
$\mathrm{Hom}(H_{\pi_{\rho_1}},H_{\pi_{\rho_2}})$.
Since $(\rho_1,\rho_2)$ is the fixed point subset of
$\cK_{\rho_2}\cK_{\rho_1}^*$  under $\theta$,
Theorem 3.5, (4) implies that the above identification carries
$(\rho_1,\rho_2)\subset (\tr_1,\tr_2)$ onto
$(\pi_{\rho_1},\pi_{\rho_2})\subset
\mathrm{Hom}(H_{\pi_{\rho_1}},H_{\pi_{\rho_2}})$.
Therefore, we define a functor $\beta$ setting $\beta_{T}$,
$T\in (\pi_{\rho_1},\pi_{\rho_2})$ to be the corresponding element in
$(\rho_1,\rho_2)$.
Note that since $\beta$ maps every irreducible object to an irreducible
object, Lemma 5.1 implies that the resulting crossed product gives
irreducible inclusion of factors $N\subset N\rtimes_\beta \hK=:M$,
and the dual action $\alpha:=\widehat{\beta}$ is minimal.

Another way to construct $(M,\alpha, K)$ from $\cD$ is to apply
the Doplicher-Roberts duality theorem \cite[Theorem 7.1]{DR1}
to $\cD$ (see also \cite{DR2}).
For this purpose, we have to specify the permutation symmetry
$\cE(\rho_1,\rho_2)$ for every pair $\rho_1,\rho_2\in \cD$.
Indeed, we set
$$\cE(\rho_1,\rho_2):=\sum_{i=1}^{d(\rho_1)}\sum_{j=1}^{d(\rho_2)}
V(\rho_2)_jV(\rho_1)_iV(\rho_2)_j^*V(\rho_1)_i^*. $$
Note that $\cE(\rho_1,\rho_2)$ belongs to $\tN^{\theta}=N$ and
it does not depend on the choice of implementing systems.
It is a routine work to show that $\cE(\rho_1,\rho_2)$ is a permutation
symmetry for $\cD$, and  that the resulting crossed product by $\cD$ and
the dual action constructed in \cite{DR1} are the same as $M$ and $\alpha$
above.

In general, $\cD$ does not uniquely determine either the compact group
$K$ or the Roberts action $\beta$ \cite{IK}.
(When $K$ is abelian, the obstruction is nothing but the usual
$H^2(\hK,\T)$ obstruction for uniqueness of liftings of the $\hK$-kernel
with a trivial 3 cocycle \cite{J1}, \cite{Su1}.)
Whenever we deal with Roberts actions consisting of modular
endomorphisms in the rest of this paper, we always assume that we make
the above choice of $\beta$ on the arrows
(or equivalently, the above choice of the permutation symmetries).

\begin{theorem} Let the notations be as above. Then,
\begin{itemize}
\item[$(1)$] $\tM$ is generated by $\tN$ and the relative commutant
$\tM\cap \tN'$,
and $\tM\cap \tN'$ coincides with the center $Z(\tM)$ of $M$.
The smooth flow of weights $F^M$ of $M$ is given by
the skew product $(X_N\times_cK,\mu_N\times m_K)$, that is:
$X_M=X_N\times K$, $\mu_M=\mu_N\times m_K$, and the flow is given by
$$(\omega,k)\cdot t=(\omega\cdot t, kc(\omega,t)),\quad t\in \R, \;k\in K,
\; \omega\in X_N.$$
The factor map from $X_M$ onto $X_N$ corresponding to the inclusion
$Z(\tN)\subset Z(\tM)$ is the projection onto the first component.
\item[$(2)$] The Connes-Takesaki module $\mo(\alpha_k)$ of $\alpha_k$, $k\in
K$
is given by
$$\mo(\alpha_k)(f)(\omega,l)=f(\omega,k^{-1}l),\quad f\in
L^\infty(X_M,\mu_M),\; (\omega,l)\in X_M.$$
\end{itemize}
\end{theorem}

\begin{proof} Before starting the proof, we simplify the notation a little.
We may and do assume $\hK\subset \cR$ thanks to the assumption (iii).
For $\rho=\beta_\pi$, $(\pi,H_\pi)\in \hK$, we simply denote
$\cK_\pi=\cK_\rho$ and $V(\pi)_i=V(\rho)_i$.

(1) We take $\cH_\pi\in M$ and an orthonormal basis
$\{W(\pi)_i\}_{i=1}^{d(\pi)}$
as before, where $\{W(\pi)_i\}_{i=1}^{d(\pi)}$ is transformed by $\alpha$ in
the
same way as the canonical basis of $H_\pi=\C^{d(\pi)}$.
Let $E$ be the unique normal conditional expectation from $M$ onto $N$, and
$\tE$ be the natural extension of $E$ to $\tM$ that is a normal conditional
expectation from $\tM$ onto $\tN$.
As in the proof of Lemma 5.1, for $x\in \tM$ we set
$$x(\pi)_i=d(\pi)\tE(W(\pi)_i^*x).$$
Then, $x$ is in $\tM\cap \tN'$ if and only if $x(\pi)_i\in (\id,\tb_\pi)$
for every $(\pi, H_\pi)\in \hK$ and every $i=1,2,\cdots, d(\pi)$.
This is further equivalent to that there exists $z(\pi)_{ij}\in Z(\tN)$
such that
$$x(\pi)_i=\sum_{j=1}^{d(\pi)}z(\pi)_{ij}V(\pi)_j,$$
for every $(\pi, H_\pi)\in \hK$ and every $i=1,2,\cdots, d(\pi)$.
Note that thanks to Proposition 4.2, (4) and Lemma 5.4, $W(\pi)_i$
commutes with every element in $Z(\tN)$.
Therefore, in the same way as in \cite[Lemma 3.8]{ILP}, we can conclude that
$\tM\cap \tN'$ is generated by $Z(\tN)$ and $\{f(\pi_{ij})\}_{\pi,i,j}$,
where $f(\pi_{ij})=W(\pi)_i^*V(\pi)_j$.
Since $\tM$ is generated by $\tN$ and $\{W(\pi)_i\}_{\pi,i}$, it is also
generated
by $\tN$ and $\{f(\pi_{ij})\}_{\pi,i,j}$ as well.
Thus, to prove $Z(\tM)=\tM\cap \tN'$, it suffices to show that
$\{f(\pi_{ij})\}_{\pi,i,j}$ generate a commutative von Neumann algebra.

Let $B$ be the von Neumann algebra generated by $\{f(\pi_{ij})\}_{\pi,i,j}$.
We show that there exists a $*$ isomorphism from
$L^\infty(K)$ onto $B$ sending $\pi_{ij}$ to $f(\pi_{ij})$.
For this, it suffices to check that both $\pi_{ij}$ and $f(\pi_{ij})$
have the same algebraic relations because we have
$\tE(f(\pi_{ij}))=m_K(\pi_{ij})
=\delta_{\pi,1}$, where
$$m_K(\xi)=\int_K\xi(k)dm_K(k), \quad \xi\in L^\infty(K).$$
First, we show $f(\pi_{ij})^*=f(\bp_{ij})$.
Let $R_\pi\in (\id,\beta_{\bp}\cdot\beta_\pi)$ and
$R_{\bp}\in (\id,\beta_\pi\cdot\beta_{\bp})$
be isometries defined by
$$R_\pi=\frac{1}{\sqrt{d(\pi)}}\sum_{i=1}^{d(\pi)}V(\bp)_iV(\pi)_i,$$
$$R_{\bp}=\frac{1}{\sqrt{d(\pi)}}\sum_{i=1}^{d(\pi)}V(\pi)_iV(\bp)_i.$$
Then, the definition of the crossed product by a Roberts action implies that
we also have
$$R_\pi=\frac{1}{\sqrt{d(\pi)}}\sum_{i=1}^{d(\pi)}W(\bp)_iW(\pi)_i,$$
$$R_{\bp}=\frac{1}{\sqrt{d(\pi)}}\sum_{i=1}^{d(\pi)}W(\pi)_iW(\bp)_i.$$
Thus, we get
\begin{eqnarray*} f(\pi_{ij})^*&=& V(\pi)_j^*W(\pi)_i
=\sqrt{d(\pi)}V(\pi)_j^*W(\bp)_i^*R_{\pi}
=\sqrt{d(\pi)}W(\bp)_i^*\tb_{\bp}(V(\pi)_j^*)R_{\pi}\\
&=&W(\bp)_i^*V(\bp)_j=f(\bp)_{ij}.
\end{eqnarray*}

For the same reason, we have
$$\cE(\pi,\sigma)=\sum_{i=1}^{d(\pi)}\sum_{a=1}^{d(\sigma)}
V(\sigma)_aV(\pi)_iV(\sigma)_a^*V(\pi)_i^*
=\sum_{i=1}^{d(\pi)}\sum_{a=1}^{d(\sigma)}
W(\sigma)_aW(\pi)_iW(\sigma)_a^*W(\pi)_i^*,$$
which implies
\begin{eqnarray*}f(\pi_{ij})f(\sigma_{ab})
&=& W(\pi)_i^*V(\pi)_jW(\sigma)_a^*V(\sigma)_b
=W(\pi)_i^*W(\sigma)_a^*\tb_\sigma(V(\pi)_j)V(\sigma)_b\\
&=&W(\pi)_i^*W(\sigma)_a^*V(\sigma)_bV(\pi)_j
=W(\pi)_i^*W(\sigma)_a^*\cE(\pi,\sigma)V(\pi)_jV(\sigma)_b\\
&=&W(\sigma)_a^*W(\pi)_i^*V(\pi)_jV(\sigma)_b=f(\sigma_{ab})f(\pi_{ij}).
\end{eqnarray*}
Thus, $B$ is abelian.
For the same reason again, there exist
${T_{(\pi,i),(\sigma,a)}^{(\nu,s)}}\in
(\beta_{\nu},\beta_\pi\cdot\beta_\sigma) $ such that
$$V(\pi)_i V(\sigma)_a=\sum_{(\nu,H_\nu)\in \hK} \sum_{s=1}^{d(\nu)}
T_{(\pi,i),(\sigma,a)}^{(\nu,s)} V(\nu)_s,$$
$$W(\pi)_i W(\sigma)_a=\sum_{(\nu,H_\nu)\in \hK} \sum_{s=1}^{d(\nu)}
T_{(\pi,i),(\sigma,a)}^{(\nu,s)} W(\nu)_s,$$
which implies
$$f(\pi_{ij})f(\sigma_{ab})=\sum_{(\nu,H_\nu)\in \hK} \sum_{s,t=1}^{d(\nu)}
\inpr{T_{(\pi,j),(\sigma,b)}^{(\nu,t)} }{T_{(\pi,i),(\sigma,a)}^{(\nu,s)} }
f(\nu_{st}).$$
Note that thanks to the Peter-Weyl theorem,
$T_{(\pi,i),(\sigma,a)}^{(\nu,s)}$ is obtained by
$$T_{(\pi,i),(\sigma,a)}^{(\nu,s)}=d(\nu)E(W(\pi)_i
W(\sigma)_aW(\nu)_s^*).$$
Thus, we have
\begin{eqnarray*}\lefteqn{
\inpr{T_{(\pi,j),(\sigma,b)}^{(\nu,t)} }{T_{(\pi,i),(\sigma,a)}^{(\nu,s)} }
}\\
&=&d(\nu)^2\int _K
E(W(\nu)_sW(\sigma)_a^* W(\pi)_i^*\alpha_g(W(\pi)_jW(\sigma)_bW(\nu)_t^*))
dm_K(g)\\
&=& d(\nu)^2\sum_{r=1}^{d(\nu)}\int _K
\pi_{ij}(g)\sigma_{ab}(g)\overline{\nu_{rt}(g)}
E(W(\nu)_sW(\nu)_r^*)dm_K(g)\\
&=&d(\nu)\int_K \pi_{ij}(g)\sigma_{ab}(g)\overline{\nu_{st}(g)}dm_K(g).
\end{eqnarray*}
This shows that
$$\pi_{ij} \sigma_{ab}=\sum_{(\nu,H_\nu)\in \hK} \sum_{s,t=1}^{d(\nu)}
\inpr{T_{(\pi,j),(\sigma,b)}^{(\nu,t)} }{T_{(\pi,i),(\sigma,a)}^{(\nu,s)}}
\nu_{st}$$
holds as well thanks to the Peter-Weyl theorem again.
Therefore, there exists a desired isomorphism from $L^\infty(K)$ onto $B$.

Let $\varphi_0$ be a faithful normal state of $\tN$ whose restriction
is given by the measure $\mu_N$, and we set $\varphi=\varphi_0\cdot \tE$.
Then, for $z\in Z(\tN)$, $f\in B$, we have
$$\varphi(zf)=\varphi_0(z\tE(f))=\varphi(z)\varphi(f).$$
This shows that $Z(\tM)=Z(\tN)\otimes B$.
$\theta$ acts on $f(\pi_{ij})$ as
$$\theta_t(f(\pi_{ij}))=W(\pi)_i^*\theta_i(V(\pi)_j)
=\sum_{l=1}^{d(\pi)}W(\pi)_i^*V(\pi)_l\pi_{lj}(c(t))
=\sum_{l=1}^{d(\pi)}f(\pi_{il})\pi_{lj}(c(t)),$$
which shows that the flow is given by the skew product.

(2) Since $\ta$ acts on $V(\pi)_i$ trivially, we have
$$\ta_k(f(\pi_{ij}))=\alpha_k(W(\pi)_i^*)V(\pi)_j
=\sum_{l=1}^{d(\pi)}\overline{\pi_{li}(k)}f(\pi_{lj})
=\sum_{l=1}^{d(\pi)}\pi_{il}(k^{-1}) f(\pi_{lj}).$$
This means that the Connes-Takesaki module of $\alpha_k$
is given as in the statement.
\end{proof}

\begin{definition} Let $N$ be a type III factor, $K$ be a compact group,
and $c\in Z^1(F^N,K)$ be a minimal cocycle with $K_c=K$.
For each $(\pi,H_\pi)\in \hK$, we choose $\rho_\pi\in \eNm$ such that
$\delta_m([\rho_\pi])=[\pi\cdot c]$.
Let $\cD$ be the subset of $\eNm$ generated by
$\{\rho_\pi\}_{(\pi,H_\pi)\in \hK}$ and satisfying the above conditions
(i)-(iv).
We denote by $N\otimes_c L^\infty(K)$ the above crossed product
$N\rtimes_\beta\hK$ and call it the \textit{skew product of $N$ and
$L^\infty(K)$ by the cocycle $c$}.
We call $\hb$ the \textit{dual action} of the skew product.
\end{definition}
Note that $(N\otimes_c L^\infty(K),\hb)$ depends only on the class
$[c]\in H^1(F^N, K)$.

Our next task is to show the ``converse" of the above theorem, and we need
some preparation.

\begin{lemma} Let $\alpha$ be an action of a locally compact group $G$ on
a von Neumann algebra $M$, and $H$ be the standard Hilbert space of $M$.
Let $\tJ$ be the modular conjugation of $M\rtimes_\alpha G
:=\pi_\alpha(M)\vee (\C\otimes l(G)'')$ acting on $L^2(G,H)$ as given by
\cite[Lemma 2.8]{Ha1}.
Then, we have
$$\tJ ((M\rtimes_\alpha G)\cap \pi_\alpha(M)')\tJ
=(M\otimes B(L^2(G)))\cap(M\rtimes_\alpha G)'.$$
\end{lemma}

\begin{proof} Let $J$ be the modular conjugation of $M$ on $H$.
Then, $\tJ$ is given by
$$(\tJ\xi)(g)=\Delta_G(g)^{-1/2}u(g)^*J\xi(g),\quad \xi\in L^2(G,H),$$
where $\Delta_G$ is the modular function of $G$ and $u(g)$ is the
canonical implementation of $G$ on $H$.
Direct computation yields
$$\tJ\pi_\alpha(x)\tJ=JxJ\otimes \C,\quad x\in M,$$
$$\tJ (u(g)\otimes r(g))\tJ=1\otimes l(g),\quad g\in G.$$
Thus, we get
\begin{eqnarray*}\tJ ((M\rtimes_\alpha G) \cap \pi_\alpha(M)')\tJ&=&
\tJ ((M\otimes B(L^2(G)))\cap (u\otimes r)(G)'\cap \pi_\alpha(M)')\tJ\\
&=&(\tJ((M'\otimes \C)\vee (u\otimes r)(G)''\vee \pi_\alpha(M))\tJ)'\\
&=&(\pi_\alpha(M) \vee (\C\otimes l(G)'')\vee (M'\otimes \C))'\\
&=&(M\rtimes_\alpha G)' \cap (M\otimes B(L^2(G))).\\
\end{eqnarray*}
\end{proof}

Let $M\supset N$ be an irreducible inclusion of factors with a normal
conditional
expectation $E$.
We assume that $\tM$ is generated by $\tN$ and $\tM\cap \tN'$.
Then, the restriction of $\tE$ to $\tM\cap \tN'$ is a conditional
expectation onto $Z(\tN)$.
For $a,b\in \tM\cap \tN'$, $\tE$ satisfies the trace property
$\tE(ab)=\tE(ba)$.
Indeed, for $x\in \tN$ in the domain of $\tau$, we have
$$\tau(\tE(ab)x)=\tau(abx)=\tau(axb)=\tau(bax)=\tau(\tE(ba)x),$$
which shows $\tE(ab)=\tE(ba)$.

\begin{lemma} Let $N$ be a type $\mathrm{III}$ factor and $A$ be a
von Neumann algebra whose center includes $Z(\tN)$.
We assume that there exist an ergodic $\R$-action $\vartheta$ on $A$
extending the flow of $Z(\tN)$, and a faithful normal conditional
expectation $\epsilon$ from $A$ onto $Z(\tN)$ such that
$\theta_t\cdot \epsilon=\epsilon\cdot \vartheta_t$, $t\in \R$
and $\epsilon(xy)=\epsilon(yx)$ for $x,y\in A$.
Then, there exists a unique irreducible inclusion $M\supset N$ of
type $\mathrm{III}$
factors with a normal conditional expectation $E$ such that
$\tM=(\tM\cap \tN')\vee \tN$ and
$$(A\supset Z(\tN), \vartheta, \epsilon)\cong
(\tM\cap \tN'\supset Z(\tN),\theta|_{\tM\cap \tN'}, \tE|_{\tM\cap \tN'}).$$
\end{lemma}

\begin{proof}
First, we show uniqueness of $M$.
Let $M$ be a von Neumann algebra including $N$ with a unique conditional
expectation $E$ from $M$ onto $N$.
We assume that $\tM$ is generated by $\tM\cap\tN'$ and $\tN$.
Then, the linear span of $(\tM\cap\tN')\tN$ is dense in $\tM$ by
assumption.
For a faithful normal state $\varphi$ on $\tN$ and $a\in \tM\cap\tN'$,
$x\in \tN$, we have
$$\varphi\cdot\tE(ax)=\varphi(\tE(a)x).$$
This means that the structure of $(\tM\supset \tN, \tE, \theta)$
is uniquely determined by $(\tN,\theta|_{\tN})$ and
$$(\tM\cap\tN'\supset Z(\tN),\theta|_{\tM\cap\tN'},\tE|_{\tM\cap\tN'}).$$
On the other hand, $(\tM\supset \tN, \tE, \theta)$ completely determine
the original inclusion  via the Takesaki duality theorem.
Therefore, such $M$ as in the statement is unique if it exists.

We set $L=\tN\otimes_{Z(\tN)} A$ as in \cite{Sau}, \cite{SZ2}.
More precisely, let $H$ be the standard Hilbert space of $\tN$.
We introduce an inner product into the algebraic tensor product
$H\odot_{Z(\tN)} A$ by
$$\inpr{\xi\odot a}{\eta\odot b}=\inpr{\epsilon(b^*a)\xi}{\eta},
\quad a,b\in A,\; \xi,\eta\in H.$$
Let $H_1$ be the completion of $H\odot_{Z(\tN)} A$, and $\Lambda$ be the
natural map from $H\odot_{Z(\tN)} A$ to $H_1$.
$\tN$ and $A$ naturally act on the first and the second tensor
components of $H_1$, and we denote by $L$ the von Neumann algebra
generated by $\tN$ and $A$ in $B(H_1)$.
In terms of direct integral, $L$ can be expressed as follows
(though we do not use it in the proof):
Let
$$\tN=\int_{X_N}^\oplus \tN(\omega)d\mu_N(\omega),$$
$$A=\int_{X_N}^\oplus A(\omega)d\mu_N(\omega),$$
be disintegration of $\tN$ and $A$ over $Z(\tN)$.
Then $L$ is given by
$$L=\int_{X_N}^\oplus \tN(\omega)\otimes A(\omega)d\mu_N(\omega).$$

Thanks to the Connes-Takesaki relative commutant theorem and Lemma 5.7,
we have $\tN\cap N'=Z(\tN)$.
We show $L\cap N'=A$ using this and applying \cite[Theorem 2.3]{SZ2}
(or alternatively, using the above direct integral expression).
To do so, we need to separate $A$ and $N$ by a type I von Neumann algebra.
Indeed, let $F$ be the commutant of $Z(\tN)$ in $B(H)$.
Then, $F$ naturally acts on $H_1$ and we have $N\subset F$ and $A\subset
F'$, where the commutant $F'$ is considered in $B(H_1)$.
Applying \cite[Theorem 2.3]{SZ2} twice, we get
\begin{eqnarray*} L\cap N'&=&(N\vee Z(\tN))' \cap (\tN\vee A)
=(N\vee Z(\tN)\vee(\tN\vee A)' )'\\
&=&(N\vee Z(\tN) \vee (F\cap \tN')\vee (F'\cap A'))'\\
&=&((F\cap (\tN\vee N')')\vee (F'\cap A'))'\\
&=&((F\cap Z(\tN)')\vee (F'\cap A'))'\\
&=&(F\vee (F'\cap A'))' \\
&=&Z(\tN)\vee A=A.
\end{eqnarray*}

We introduce a 1-parameter automorphism group $\Theta$ on $L$ extending
both $\theta$ and $\vartheta$.
Let $u(t)$ be the canonical implementation of $\theta_t$ on $H$.
We set
$$U(t)\Lambda(\xi\odot a)=\Lambda(u(t)\xi \odot \vartheta_t(a)),
\quad \xi\in H,\; a\in A.$$
Then, $U(t)$ extends to a 1-parameter unitary group on $H_1$ satisfying
$$U(t)xaU(t)^*=\theta_t(x)\vartheta(a),\quad x\in \tN,\; a\in A.$$
We define $\Theta_t$ to be the restriction of $\Ad(U(t))$ to $L$.

Next, we show that there exists a faithful normal conditional expectation
$\tep$ from $L$ onto $\tN$ satisfying
$$\tep(a)=\epsilon(a),\quad a\in A.$$
Let $\varphi$ be a faithful normal state on $\tN$,
and $\xi_\varphi\in H$ be the cyclic and separating vector of
$\tN$ giving $\varphi$.
We set $\Omega=\Lambda(\xi_\varphi\odot 1)$, and claim that
$\Omega$ is cyclic and separating vector for $L$.
Indeed, $\Omega$ is clearly a cyclic vector.
To show that it is separating, we introduce a conjugation $J$ on
$H_1$ by
$$J\Lambda(\xi\odot a)=\Lambda(J_{\tN}\xi\odot a^*),
\quad \xi\in H,\; a\in A\,$$
where $J_{\tN}$ is the modular conjugation of $\tN$.
It is easy to show that $J$ extends to a conjugation such that
$JLJ\subset L'$.
Since $JLJ$ is cyclic for $\Omega$, $\Omega$ is a separating for $L$.
We denote by $\omega$ the state of $L$ defined by
$\omega(x)=\inpr{x\Omega}{\Omega}$, $x\in L$.
Then, it is a routine work to show that $J$ is the modular conjugation for
$\omega$, and the modular automorphism group
$\{\sigma^{\omega}_t\}_{t\in \mathbb{R}}$ is given by
$$\sigma^\omega_t(xa)=\sigma^\varphi_t(x)a,\quad x\in \tN,\; a\in A.$$
Thanks to Takesaki's theorem \cite{T1}, there exists a $\omega$-preserving
normal conditional expectation $\tep$ from $L$ onto $\tN$.
Since $\omega(xa)=\omega(x\epsilon(a))$ for $x\in \tN$, $a\in A$, $\tep$
satisfies
$\tep(a)=\epsilon(a)$.

We set $M=L^\Theta$ and $E$ to be the restriction of $\tep$ to $M$.
Since $\tep$ satisfies $\tep\cdot\Theta=\theta\cdot\tep$, $E$ is a faithful
normal conditional expectation from $M$ onto $N$.
Therefore, $M$ is of type III \cite{St}.
Let $\nu$ be the restriction of $\Ad(\lambda^\psi(t))$ to $M$, where
$\psi$ is a faithful normal semifinite weight on $N$.
Then, thanks to the Landstad theorem \cite[Chapter II, Section 2]{NT},
$L=M\rtimes_\nu \R$ and $\Theta$ is the dual action of $\nu$.
We claim that $M\cap N'=\C$, and in particular $M$ is a factor.
Indeed, we have
$$M\cap N'\subset L\cap N'=A.$$
Since $\Theta$ acts on $A$ ergodically, we get $M\cap N'=\C$.

Instead of showing that $\nu$ comes from the modular automorphism group
directly, we show it for the dual action of $\Theta$.
We define $\tau_1=\tau\cdot\tep$, which is a faithful normal semifinite
trace on $L$ satisfying $\tau_1\cdot\Theta_t=e^{-t}\tau_1$.
Indeed, since $\tau$ is a trace, $\sigma^{\tau_1}_t$ acts on $\tN$
trivially.
On the other hand, since $\tep$ is a faithful normal conditional
expectation from $L$ onto $\tN$ and $L\cap \tN'=A$, the restriction of
$\sigma^{\tau_1}_t$ to $A$ is the same as the restriction of
$\sigma^{\omega}_t$ to $A$ \cite{CD}, which is trivial as we saw before.
Thus, $\tau_1$ is a trace with the scaling property
$$\tau_1\cdot \Theta_t=\tau\cdot \theta_t\cdot \tep=e^{-t}\tau_1.$$
This means that the dual action of $\Theta$ is the modular automorphism
group for the dual weight $\widehat{\tau}_1$.
Thus, identifying $L\rtimes_\Theta\R$ with $M\otimes B(L^2(R))$ using
the Takesaki duality theorem, we get
$$\widehat{\tau}_1=\widehat{\tau}\cdot (E\otimes id),$$
$$\sigma^{\widehat{\tau_1}}_t=(\nu_t\otimes \Ad(r(t))).$$
Let $h$ be the generator of $r(t)$.
Then, we have $\widehat{\tau}=\psi \otimes \mathrm{Tr}(h\cdot)$, which
implies $\widehat{\tau}_1=\psi\cdot E\otimes \mathrm{Tr}(h\cdot)$.
Thus, we get $\nu_t=\sigma^{\psi\cdot E}_t$,
which implies $L=\tM$ and $\tep=\tE$.
\end{proof}

\begin{theorem}[Automatic Minimality]
Let $\alpha$ be an action of a compact group $K$ on a
type $\mathrm{III}$ factor $M$.
If the kernel of $\mo\cdot\alpha$ is trivial, then $\alpha$ is minimal.
The smooth flow of weight of the fixed point algebra $M^\alpha$ is the
factor flow of $F^M$ by $\mo\cdot\alpha_K$.
There exists a minimal cocycle $c\in Z^1(F^{M^\alpha},K)$ with
$K_c=K$, unique up to equivalence, such that
$M$ is the skew product $M^\alpha\otimes_c L^\infty(K)$, and $\alpha$ is the
dual action of the skew product.
\end{theorem}

\begin{proof} We use the same notation as in Lemma 5.4.
Lemma A.1 implies that the smooth flow of weights of $M$ is of the form
$(X_M,\mu_M)=(Y\times K,\nu\times m_K)$ with the flow and $\mo(\alpha_k)$
given by
$$(\omega,l)\cdot t=(\omega\cdot t,lc(\omega,t)),\quad (\omega,l)\in X_M,
\;t\in \R,$$
$$(\mo(\alpha_k)f)(\omega,l)=(\omega,k^{-1}l),\quad k\in K,\;f\in
L^\infty(X_M,\mu_M),$$
where $c$ is a minimal cocycle with $K_c=K$.
According to the above splitting of $(X_M,\mu_X)$,
we set $C=L^\infty(K)\subset Z(\tM)$ and $Z=L^\infty(Y,\nu)\subset Z(\tM)$
that is the fixed point subalgebra of $Z(\tM)$ under $\ta$.
Note that the factor flow $(Z,\theta|_Z)$ and the class of $c$ are uniquely
determined by $F^M$ and $\mo\cdot\alpha$.
Since $C$ is an equivariant copy of $L^\infty(K)$,
the Landstad's type theorem for coactions \cite[Chapter II, Section 2]{NT}
implies that $\tM$ is the crossed product $\tN\rtimes \hK$ by a coaction
of $K$, where we use the fact that $\tN$ is the fixed point subalgebra of
$\tM$ under $\ta$.
Moreover, since $\tN$ commutes with $C$, the coaction is trivial and we
actually have $\tM\cong\tN\otimes C$, where $\tN$ and $C$ are identified
with $\tN\otimes \C$ and $\C\otimes C$ in the right-hand side respectively.
This implies $Z=Z(\tN)$, $\tM\cap \tN'=Z(\tM)$ and
$\tM=(\tM\cap \tN')\vee \tN$.
Since $\theta$ acts on $Z$ ergodically, $N$ is a factor.
If $N$ were semifinite, $c$ would be a coboundary, which would imply that
$K$
is trivial and $M=N$.
This is contradiction because $M$ is of type III.
Therefore, $N$ is a type III factor.
Now, Theorem 5.5 and Lemma 5.8 show that $M$ is the skew product
$N\otimes_c L^\infty(K)$.
Since the Galois group for $M\supset N$ coincides with $\alpha_K$
\cite{R},\cite{AHKT} and every member $\alpha_k$ in the Galois group
is specified by $\mo(\alpha_k)$,
we conclude that $\alpha$ is the dual action of the skew product.
\end{proof}

\begin{corollary} If $M$ is an AFD factor, then for every compact subgroup
$K$ of $\mathrm{Aut}(F^M)$, there exists a minimal action $\alpha$ of $K$ on
$M$
such that $\mo\cdot \alpha_k=k$ for all $k\in K$.
$\alpha$ is unique up to conjugacy: namely if $\alpha^1$ also satisfies
the condition, then there exists an automorphism $\nu\in \mathrm{Aut}(M)$
with a trivial Connes-Takesaki module such that
$$\alpha_k=\nu\cdot\alpha^1_k\cdot \nu^{-1},\quad k\in K.$$
\end{corollary}

\begin{proof} This is an easy consequence of Theorem 5.5, Theorem 5.9, and
the fact that AFD type III factors are determined by their flows
\cite{C2}, \cite{C3}, \cite{C4}, \cite{Kr}, \cite{Ha3}.
\end{proof}

As for the existence of a homomorphic lifting in the above corollary,
a more general statement is known: there always exists a homomorphic lifting
$\alpha$ from $\mathrm{Aut}(F^M)$ to $\mathrm{Aut}(M)$ for every AFD
type III factor $M$ \cite{ST2}.

Now we characterize (not necessarily finite index) subfactors
whose canonical endomorphism is decomposed into modular endomorphisms.
Let $M\supset N$ be an irreducible inclusion of factors with a unique normal
conditional expectation from $M$ onto $N$.
We denote by $M_1$ and $\hE$ the basic extension and the dual operator
valued weight from $M_1$ to $M$ respectively \cite{Ha2}, \cite{K1},
\cite[Section 2]{ILP}.
We say that the inclusion $M\supset N$ is \textit{discrete} if the
restriction of $E\cdot \hE$ on the relative commutant $M_1\cap N'$
is semifinite.
We say that a discrete inclusion is \textit{unimodular} if
$E\cdot \hE|_{M_1\cap N'}$ is invariant under $j(x):=J_Mx^*J_M$,
where $J_M$ is the modular conjugation of $M$
(see \cite[Section 3]{ILP} for conditions equivalent to
these.)

\begin{theorem} Let $M\supset N$ be an inclusion of type $\mathrm{III}$
factors with a unique normal conditional expectation $E$ from $M$ onto $N$.
Then, the following conditions are equivalent:
\begin{itemize}
\item[$(1)$] $M\supset N$ is a unimodular inclusion and the restriction
of the canonical endomorphism $\gamma:M\longrightarrow N$ to
$N$ is decomposed into modular endomorphisms.
\item[$(2)$] $M\supset N$ is a discrete inclusion and $\tM$ is generated by
$\tN$ and $\tM\cap \tN'$.
\item[$(3)$] $(\tM\cap \tN',\theta)$ is an extension of the flow
$(Z(\tN),\theta)$ with relatively discrete spectrum in the sense of
Definition A.2, and $\tM$ is generated by $\tN$ and $\tM\cap \tN'$.
\item[$(4)$] There exists a compact group $K$, a faithful ergodic action
$\Psi$ of
$K$ on a von Neumann algebra $B$, a minimal action $\alpha$ of $K$ on a
factor $P$
such that the kernel of $\mo\cdot \alpha$ is trivial and
$$(M\supset N)\cong ((P\otimes B)^{\alpha\otimes \Psi}\supset
P^\alpha\otimes \C).$$
\end{itemize}

In the case (4), $K$ is determined by the cocycle $c$ corresponding to the
system
of modular endomorphisms appearing in the decomposition of $\gamma|_N$,
$P$ is the skew product $N\otimes_c L^\infty(K)$, and $\alpha$ is the dual
action
of the skew product.
$B$ and $\Psi$ are uniquely determined from
$(\tM\cap \tN'\supset Z(\tN),\theta)$ by Theorem A.4.

When one (and all) of the above conditions holds, $N$ is AFD
if and only if $M$ is AFD.
\end{theorem}

\begin{proof} (1) $\Rightarrow$ (2). We assume that (1) holds and the
irreducible decomposition of $\gamma|_N$ is given by
$$[\gamma|_N]=\bigoplus_{\xi\in \Xi}n_\xi[\rho_\xi].$$
Thanks to \cite[pp. 39]{ILP}, the multiplicity $n_\xi$ of $\rho_\xi$
is always finite.
We set
$$\cL_\xi=\{X\in M;\; Xx=\rho_\xi(x)X,\; \forall x \in N\}.$$
Then, the dimension of $\cL_\xi$ is $n_\xi$ and $N$ is generated by $N$
and $\bigcup_{\xi\in \Xi} \cL_\xi$ \cite[Theorem 3.3]{ILP}.
Let $\{X(\xi)_i\}_{i=1}^{n_\xi}$ be an orthonormal basis of $\cL_\xi$.
We choose a dominant weight $\psi_0$ on $N$ and assume that
$(\psi_0,\rho_\xi)$
is an invariant pair for every $\xi\in \Xi$ as before.
Then, since $M\supset N$ is unimodular,
\cite[Theorem 3.3, Remark 3.4]{ILP} imply that $\cL_\xi$ is in the
centralizer of $\psi=\psi_0\cdot E$.
By assumption, there exists an implementing system
$\{V(\xi)_i\}_{i=1}^{d(\xi)}\subset \tN$ for $\tr_\xi$.
As $X(\xi)_i$ is in the centralizer of $\psi$, it commutes with
$\lambda^\psi(t)$, and so
$$V(\xi)_i^*X(\xi)_j\in \tM\cap \tN',\quad \xi\in \Xi,
\;1\leq i\leq d(\rho_\xi),\; 1\leq j\leq n_\xi. $$
Since $\tM$ is generated $\tN$ and $\cup_{\xi\in \Xi} \cL_\xi$, it is also
generated by $\tN$ and $\tM\cap\tN'$.

(2) $\Rightarrow$ (3). We assume (2), and set $A=\tM\cap \tN'$,
$Z=Z(\tN)$, $\epsilon =\tE|_A$, and $\vartheta_t=\theta_t|_A$.
Then, $\epsilon$ is a normal conditional expectation with the
trace property $\epsilon(ab)=\epsilon(ba)$, $a,b\in A$.
Let $M_1$ be the basic extension of $M\supset N$.
As before, we regard $\tM$ as a subalgebra of $\tM_1$ identifying
$\lambda^\psi(t)$ with $\lambda^{\psi\cdot \hE}(t)$.
Then, $\tM_1$ is the basic extension of $\tM\supset \tN$ with respect to
$\tE$ \cite[Lemma 2.4]{ILP}, \cite[Theorem 2.1]{FI}.

On the other hand, we consider the inclusion $Z\subset A$ acting on the
standard Hilbert space $L^2(A)$ of $A$.
Then, the basic extension of this inclusion is $Z_1=Z'\subset B(L^2(A))$.
Let $u$ be the canonical implementation of $\vartheta$, and $J_A$ be
the modular conjugation of $A$.
We use the same notation as in the proof of Lemma 5.8, and identify
$\tM$ with $L$.
Note that $Z_1$ naturally acts on $H_1$ through the action on the second
tensor component.
In the same way as before, we get the following:
$$\tM_1=J\tN'J=J_{\tN} F J_{\tN}\vee J_AZ_1J_A=\tN\vee Z_1,$$
where everything takes place in $B(H_1)$.
Thus, $\tM_1\cap \tN'=Z_1,$ and
$$(\tM_1\cap \tN', \theta,E\cdot\hE)=
(Z_1, \Ad(u(t)),\epsilon \cdot \widehat{\epsilon}),$$
$$((M_1\cap N')_{E\cdot\hE},E\cdot\hE)=
(Z_1\cap u(\R)',\epsilon \cdot \widehat{\epsilon}),$$
where $(M_1\cap N')_{E\cdot\hE}$ is the fixed point subalgebra of
$M_1\cap N'$ under the relative modular automorphism group
$\{\sigma^{E\cdot \hE}_t\}_{t\in \R}$ \cite{Ha2}.
Since $M\supset N$ is discrete, \cite[Theorem 3.3]{ILP} implies that
$(M_1\cap N')_{E\cdot\hE}$ is a direct sum of matrix algebras and the
restriction of $E\cdot\hE$ to $(M_1\cap N')_{E\cdot\hE}$ is still
semifinite.
Therefore, if $\epsilon \cdot \widehat{\epsilon}$ is the center valued trace
$\natural$ of $Z_1$, $(A,\vartheta)$ is an extension of the
flow $F^N$ with relatively discrete spectrum.
Indeed, let $\mu$ be a faithful normal state on $Z(\tN)$.
Then using the spatial derivative, we can show that the modular automorphism
group of $\mu\cdot \epsilon \cdot \widehat{\epsilon}$ is implemented by
the modular operator $\Delta_{\mu\cdot \epsilon}^{it}$ of $\mu\cdot
\epsilon$ on $L^2(A)$ \cite{K1}.
However, since $\mu\cdot \epsilon$ is a trace and $\Delta_{\mu\cdot
\epsilon}$
is trivial, $\mu\cdot \epsilon \cdot \widehat{\epsilon}$ is a trace as well.
This means that $\epsilon \cdot \widehat{\epsilon}=\natural(z\cdot )$
holds with some density $z$ affiliated with $Z(\tN)$.
Let $e_{Z(\tN)}$ be the Jones projection for $\epsilon$, and $a\in A$.
Then, $ae_Aa^*$ is considered as a field of scalar multiples of rank one
projections in the central decomposition of $Z_1$.
Since
$$\epsilon \cdot \widehat{\epsilon}(ae_Aa^*)
=\epsilon(aa^*)=\epsilon(a^*a),$$
we actually have $z=1$.

(3) $\Rightarrow$ (4). We assume (3).
Thanks to Theorem A.4, there exists a compact group $K$, a minimal
cocycle $c\in Z^1(F^N,K)$ with $K_c=K$, and a faithful ergodic action
$\Psi$ of $K$ on a von Neumann algebra $B$ such that
$$(\tM\cap\tN'\supset Z(\tN), \theta)\cong
((Z(\tN)\otimes B)\supset (Z(\tN)\otimes \C),\Psi^c),$$
(see Theorem A.4 for the notations).
Let $P$ be the skew product $N\otimes_c L^\infty(K)$ and $\alpha$ be
the dual action.
We set $Q=(P\otimes B)^{\alpha\otimes \Psi}$.
Then, thanks to ergodicity of $\Psi$, $Q\supset N$ is an irreducible
inclusion of factors.
Let $\tau_B$ be the unique $\Psi$-invariant trace \cite{HLS}.
Then, the restriction of $(\id\otimes \tau_B)$ to $Q$ gives a normal
conditional expectation from $Q$ onto $N$.
Thanks to Lemma 5.8, to prove the statement it suffices to show that
$\tQ$ is generated by $\tN$ and $\tQ\cap \tN'$, and
$$(\tQ\cap \tN'\supset Z(\tN), \theta)\cong
(Z(\tN)\otimes B\supset Z(\tN)\otimes \C,\Psi^c).$$
Lemma 5.4 and Theorem 5.9 show that $\tQ$ is actually given by
$$\tQ=(\tP\otimes B)^{\ta\otimes \Psi}
=\tN\otimes (L^\infty(K)\otimes B)^{Ad(l)\otimes \Psi},$$
where $l$ is the left regular representation.
Thus, we have
$$\tQ\cap \tN'=L^\infty(X_N,\mu_N)\otimes
(L^\infty(K)\otimes B)^{\Ad(l)\otimes \Psi},$$
which shows $\tQ=\tN\vee(\tQ\cap \tN')$.
Note that $\theta$ acts on $B$ trivially and it acts on
$Z(\tP)=L^\infty(X_N,\mu_N)\otimes L^\infty(K)$ as in Theorem
5.9.

Now we assume that $Z(\tP)\otimes B$ acts on
$L^2((X_P,\mu_P), L^2(B))$,
where $L^2(B)$ is the standard Hilbert space of $B$ and
$(X_P,\mu_P)=(X_N\times K,\mu_N\times m_K)$.
Let $u$ and $v$ be the canonical implementations of $F^{P}$ on
$L^2(X_P,\mu_P)$ and $\Psi$ on $L^2(B)$ respectively.
Then, $u$ is given by
$$(u(t)\xi)(\omega,k)=\frac{dt\cdot \mu_P}{d\mu_P}(\omega)^{1/2}
\xi(\omega\cdot t,kc(\omega,t)), \quad \xi\in L^2(X_P,\mu_P).$$
We introduce a unitary operator $W$ on $L^2(K,L^2(B))$ by
$$(W\eta)(k)=v(k)\eta(k),\quad \eta\in L^2(K,L^2(B)),$$
which satisfies $W(l(k)\otimes 1)W^*=l(k)\otimes v(k).$
Let $\Pi$ be the restriction of $\Ad(W)$ on $L^\infty(K)\otimes B$,
which is an automorphism of $L^\infty(K)\otimes B$.
Then, the above relation implies
$$(L^\infty(K)\otimes B)^{\Ad(l)\otimes \Psi}=\Pi(\C\otimes B),$$
$$(r(k)\otimes 1)\Pi(1\otimes b)(r(k)^*\otimes 1)
=\Pi(1\otimes \Psi_k(b)),\quad b\in B.$$
For $f\in L^\infty(X_N,\mu_N)$, $b\in B$, and
$\xi\in L^2((X_N,\mu_N), L^2(K)\otimes L^2(B))$ we have
\begin{eqnarray*} \lefteqn{
(u(t)(f\otimes \Pi(1\otimes b))u(t)^*) \xi(\omega)}\\
&=&f(\omega\cdot t)
(r(c(\omega,t))\otimes 1)\Pi(1\otimes b)
(r(c(\omega,t))^*\otimes 1)\xi(\omega)\\
&=&f(\omega\cdot t) \Pi(1\otimes \Psi_{c(\omega,t)}(b))\xi(\omega).
\end{eqnarray*}
Thus, we get
$$(\tQ\cap \tN'\supset Z(\tN), \theta)\cong
(Z(\tN)\otimes B\supset Z(\tN)\otimes \C,\Psi^c).$$

(4) $\Rightarrow$ (1). We assume (4) and set $N=P^\alpha$,
$M=(P\otimes B)^{\alpha\otimes \Phi}$.
Let $\tau_B$ be the unique $\Psi$-invariant trace on $B$, $L^2(B)$ be the
GNS Hilbert space of $\tau_B$, and $\Lambda_{\tau_B}$ be the natural map
from $B$ into $L^2(B)$.
We denote by $v$ the canonical implementation of $\Psi$ on $L^2(B)$.
For each $(\pi,H_\pi)\in \hK$, we denote by $B_\pi$ the spectral subspace
corresponding to $\pi$, which is the image of the map
$$E_\pi(x)=d(\pi)\int_{K}\mathrm{Tr}(\overline{\pi}(k))\Psi_k(x)dm_K(k),
\quad x\in B.$$
It is known that the multiplicity of $B_\pi$ is finite \cite{HLS} and
there exists an orthonormal basis $\{b(\pi,s)_i\}\subset B_\pi$,
$1\leq s\leq n_\pi$, $1\leq i\leq d(\pi)$ with respect to the inner product
of $L^2(B)$ such that
$$\Psi_k(b(\pi,s)_i)=\sum_{j=1}^{d(\pi)}\pi_{ji}(k)b(\pi,s)_j,\quad k\in
K.$$
Since $\tau_B$ is given by the average of $\Psi_k$ over $K$, the Peter-Weyl
theorem implies
$$\sum_{i=1}^{d(\pi)}b(\pi,s)_i^*b(\pi,t)_i
=\sum_{i=1}^{d(\pi)}b(\pi,t)_ib(\pi,s)_i^*=d(\pi)\delta_{s,t}.$$
Note that $\{b(\pi,s)_i\}_{\pi,s,i}$ form an orthonormal basis of $L^2(B)$.
Though expansion of $x\in B$ with this basis does not converge in weak
topology in general, we can still show that $\{b(\pi,s)_i\}_{\pi,s,i}$
generate $B$ using the same argument as in \cite[Lemma 3.8]{ILP}.
The point of the argument is that when a von Neumann subalgebra is globally
invariant under the modular automorphism group of a state, the corresponding
subspace in the GNS subspace determines the subalgebra.
We use this principle in the rest of the proof without mentioning it.

Let $E_0$ be the unique normal conditional expectation from $P$ onto
$N$, $\varphi_0$ be a faithful normal state on $N$, and
$\varphi=\varphi_0\cdot E_0$.
We denote by $L^2(P)$, $\Lambda_\varphi$, and $w$ the GNS Hilbert space of
$\varphi$, the natural map from $P$ into $L^2(P)$, and the canonical
implementation of $\alpha$ on $L^2(P)$ respectively.
We regard the subspace of $w\otimes v$-invariant vectors in
$L^2(P)\otimes L^2(B)$ as the GNS Hilbert space $L^2(M)$ of the restriction
of $\varphi\otimes \tau_B$ to $M$.
Let $\beta$ be the Roberts action of $\hK$ on $N$ whose dual action is
$\alpha$, and $\{W(\pi)_i\}_{i=1}^{d(\pi)}\subset P$, $(\pi,H_\pi)\in \hK$
be an orthonormal basis of $\cH_\pi\subset P$ as before.
Thanks to Theorem 5.9, $\beta_\pi$ is a modular endomorphism for all
$(\pi,H_\pi)\in \hK$.
Note that $N \cH_\pi$ is globally invariant under the modular
automorphism group $\{\sigma^\varphi_t\}_{t\in \R}$
\cite[Lemma 3.8]{ILP}.
Let
$$\cL_\pi=\{X\in M;\; Xx=\beta_\pi(x)X,\; \forall x\in N \}.$$
Thanks to \cite[Theorem 3.3]{ILP}, if $M$ is generated by $N$ and
$\bigcup_{(\pi,H_\pi)\in \hK} \cL_\pi$, $M\supset N$ is a discrete
inclusion.

Now we set
$$X(\pi)_s=\frac{1}{\sqrt{d(\pi)}}
\sum_{i=1}^{d(\pi)}W(\pi)_i\otimes b(\pi,s)_i^*\in \cL_\pi.$$

Then, $\{X(\pi)_s\}_{s=1}^{n_\pi}$ are isometries with mutually orthogonal
ranges.
The definition of the crossed product by a Roberts action implies
$$L^2(P)=\bigoplus_{(\pi,H_\pi)\in \hK}
\overline{\Lambda_\varphi(\cH_\pi^* N)}^{||\cdot||}.$$
Therefore, direct computation using the orthonormal basis
$\{b(\pi,s)_i\}_{\pi,s,i}$ yields
$$L^2(M)=\bigoplus_{(\pi,H_\pi)\in \hK}
\overline{\Lambda_{\varphi\otimes \tau_B}(\cL_\pi^*N)}^{||\cdot||}.$$
Thus, $M$ is generated by $N$ and $\bigcup_{(\pi,H_\pi)\in \hK} \cL_\pi$,
and so $M\supset N$ is a discrete inclusion.

Let $E$ be the restriction of $\id\otimes \tau_B$ to $M$, which is a
conditional
expectation from $M$ onto $N$.
Then, we have
\begin{eqnarray*}
d(\pi)E(X(\pi)_sX(\pi)_t^*)&=&\sum_{i,j=1}^{d(\pi)}
\tau_B(b(\pi,s)_i^* b(\pi,t)_j)W(\pi)_iW(\pi)_j^*\\
&=&\delta_{s,t}\sum_{i=1}^{d(\pi)}W(\pi)_iW(\pi)_i^*\\
&=&\inpr{X(\pi)_s}{X(\pi)_t}.
\end{eqnarray*}
This shows that the inclusion is unimodular \cite[Theorem 3.3]{ILP}.

When $N$ is AFD, $\tN$ is AFD and $\tM$ is AFD thanks to the
proof of Lemma 5.8 and Theorem A.4.
Therefore, $M$ is AFD.
Since there exists a conditional expectation from $M$ onto
$N$, the converse also holds \cite{C3}.
\end{proof}

\begin{remark} In the above, if $\tM\cap \tN'$ is commutative, $B$ in (4) is
commutative, and there exists a unique closed subgroup $L\subset K$ such
that
$B=L^\infty(K/L)$ where $K$ action is given by the right translation.
Thus,
$$(M\supset N)\cong (P^{\alpha|_L}\supset P^\alpha).$$
This characterizes an inclusion coming from a (commutative) ergodic
extension of the flow of weights.
For related topics, see \cite{HK1}, \cite{HK2}, \cite{I3}, \cite{HK3},
\cite{K2}, \cite{KS}, \cite{Su2}, \cite{Sch}, \cite{Was2}.
Under the condition that the kernel of $\mo\cdot \alpha$ is trivial,
(and $L$ does not contain non-trivial normal subgroups, which is regarded as
a trivial constraint), $L\subset K$ are unique up to conjugacy.
However, if we drop this condition, there are essentially different
description of $M\supset N$ in terms of a group-subgroup subfactor
\cite{KSu}, \cite{IK}.

Another easy but interesting case is when $B$ is a matrix algebra, that
is,  $M\supset N$ is a Wassermann's type subfactor.
In this case, $M$ and $N$ have the common flow and graph change occurs
(c.f. \cite{Mas}).
\end{remark}

\subsection{Connes-Takesaki modules for Roberts actions}
In this subsection, we show how the Connes-Takesaki modules of endomorphisms
give constraints to minimal actions of compact groups, in particular
semisimple compact Lie groups.

Let $K$ be a compact group and $\beta$ be a Roberts action
of $\hK$ on a factor $N$ such that $\beta_\pi$ is irreducible and has a
Connes-Takesaki module for every $(\pi,H_\pi)\in \hK$.
We denote $m_\beta(\pi)=\mo(\beta_\pi)$ for simplicity.
Then, thanks to Proposition 4.2, $\{m_\beta(\pi)\}_{(\pi,H_\pi)\in \hK}$
form a group, which is denoted by $m_\beta(\hK)$.
We treat $m_\beta(\hK)$ as a discrete group, and we denote by
$\widehat{m_\beta}$ its dual group.

\begin{theorem} Let $K$ be a compact group and $\beta$ be a Roberts action
of $\hK$ on a factor $N$ such that $\beta_\pi$ is irreducible for every
$(\pi,H_\pi)\in \hK$.
Then,
\begin{itemize}
\item [$(1)$] If $K$ is connected, $\beta_\pi$ has a Connes-Takesaki module
for every $(\pi,H_\pi)\in \hK$.
\item [$(2)$] If $\beta_\pi$ has a Connes-Takesaki module for every
$(\pi,H_\pi)\in \hK$, then there exists a continuous isomorphism $\nu$
from $\widehat{m_\beta}$ into the center $Z(K)$ of $K$ satisfying
$$\chi(m_\beta(\pi))1_\pi=\pi(\nu(\chi)),\quad \chi\in \widehat{m_\beta},
\; (\pi,H_\pi)\in \hK.$$
\end{itemize}
\end{theorem}

\begin{proof} (1). Suppose that there exists $(\pi,H_\pi)\in \hK$ such that
$\beta_\pi$ does not have a Connes-Takesaki module.
Let $N=\beta_\pi(M)$, and we consider the Loi-Hamachi-Kosaki decomposition
of $M\supset N$ \cite{Loi}, \cite{HK1}, \cite{HK2} as follows:
$$M\supset P\supset Q\supset N,$$
where $P$ and $Q$ are determined by the conditions
$$\tQ=Z(\tM\cap \tN')\vee \tN,$$
$$\tP=\tM\cap Z(\tM\cap \tN')'.$$
Since $\beta_\pi$ does not have a Connes-Takesaki module, either
$M\neq P$ or $Q\neq N$ occurs.
We may assume $M\neq P$ by considering $\beta_{\overline{\pi}}$ if
necessary.
Then, $\beta_\pi \cdot \beta_{\overline{\pi}}$ contains the
the canonical endomorphism $\gamma$ for $M\supset P$.
Note that the dual inclusion $P_1\supset M$ satisfies the condition of
Theorem 5.11 and $(Z(\tP_1\cap \tM')=Z(\tP_1),\theta)$ is
a finite to one extension of the flow $F^M$ \cite{HK1}, \cite{HK2}.
Thus, Theorem 5.11 implies that $\gamma$ generates a Roberts action of
the dual $\widehat{G}$ of a finite group $G$, which is a subset of
$\{\beta_{\pi}\}_{(\pi,H_\pi)\in \hK}$.
We regard $\widehat{G}$ as a subset of $\hK$.
We define a closed normal subgroup $L\subset K$ by
$$L=\bigcap_{(\pi, H_\pi)\in \widehat{G}}\mathrm{Ker}(\pi).$$
Since $\widehat{G}$ is a finite set of irreducible representations
closed under complex conjugate and irreducible
decomposition of products of any members in $\widehat{G}$,
the quotient group $K/L$ is finite \cite[Section 30]{HR}.
However, this is contradiction because $K$ is connected.
Therefore, $\beta_\pi$ has a Connes-Takesaki module for every
$(\pi,H_\pi)\in \hK$.

(2). For $\chi\in \widehat{m_\beta}$, we consider a family of unitary
operators
$$\{\chi(m_\beta(\pi))1_\pi\}_{(\pi,H_\pi)\in \hK}.$$
Then, thanks to Theorem 4.2 and the Tannaka duality theorem
\cite[Section 30]{HR}, there exists a unique element $\nu(\chi)\in K$
such that $\chi(m_\beta(\pi))1_\pi=\pi(\nu(\chi))$ holds.
Since $\pi(\nu(\chi))$ is a scalar for every irreducible $\pi$,
$\nu(\chi)$ is in the center $Z(K)$.
$\nu$ is clearly a homomorphism from $\widehat{m_\beta}$ into $Z(K)$,
which is injective as the set of irreducible representations separates
points of $K$.
$\nu$ is continuous because the topology of $K$ is the same as that
induced by the weak topology of the image of $K$ by the direct sum
representation of all members in $\hK$.
\end{proof}

\begin{corollary} Let $\alpha$ be a minimal action of a compact connected
semisimple Lie group $K$ on a type $\mathrm{III}$ factor $M$.
Then,
\begin{itemize}
\item [$(1)$] If $M$ is of type $\mathrm{III}_1$, so is $M^{\alpha}$.
\item [$(2)$] If $M$ is of type $\mathrm{III}_\lambda$, $0<\lambda<1$,
there exists a positive integer $n$ such that $M^\alpha$ is of type
$\mathrm{III}_{\lambda^n}$.
Let $\beta$ be the Roberts action of $\hK$ on $M^\alpha$ whose dual action
is $\alpha$, then $m_\beta(\hK)\cong \widehat{m_\beta}\cong \Z/n\Z$.
\item [$(3)$] If $M$ is of type $\mathrm{III}_0$, so is $M^{\alpha}$.
(The flow of $M$ and $M^\alpha$ could be very much different in this case.)
\end{itemize}
\end{corollary}

\begin{proof} Since the center of a semisimple Lie group is a finite group,
Proposition 5.2, (2) implies that $N:=M^\alpha$ is of type III.
As before, we regard $\tN$ as a subalgebra of $\tM$ using a unique
conditional
expectation $E$.
Let $\beta$ be the Roberts action of $\hK$ on $N=M^\alpha$
whose dual action is $\alpha$.

We first assume that $N$ is of type $\mathrm{III}_\lambda$,
$0<\lambda\leq 1$.
Since every non-trivial modular endomorphism of $N$ is a composition
of an inner automorphism and a modular automorphism in this case, no
non-trivial
modular endomorphisms appear in $\{\beta_\pi\}_{(\pi,H_\pi)\in \hK}$
because a connected semisimple Lie group has no non-trivial 1 dimensional
representation.
Thus, the same computation as in the proof of Theorem 5.5 implies
$$\tM\cap \tN'=Z(\tM).$$
Thanks to Theorem 5.13, $\beta_\pi$ has a Connes-Takesaki module for every
$(\pi, H_\pi)\in \hK$, which implies
$$Z(\tM)=Z(\tN)^{m_\beta(\hK)}.$$
This shows that if $N$ is of type $\mathrm{III}_1$, so is $M$.
When $N$ is of type $\mathrm{III}_\lambda$, $0<\lambda<1$, we have
$\mathrm{Aut}(F^N)\cong \T$.
Since $Z(K)$ is a finite group, $\widehat{m_\beta}$ is finite and its dual
group $m_\beta(\hK)$ is a finite subgroup of $\T$.
Thus, $m_\beta(\hK)$ is a finite cyclic group and so is $\widehat{m_\beta}$.

Now we assume that $N$ is of type $\mathrm{III}_0$.
Then, in the same way as above we can show
$$Z(\tM)\supset Z(\tN)^{m_\beta(\hK)}.$$
Since $m_\beta(\hK)$ is a finite group, $(Z(\tN)^{m_\beta(\hK)},\theta)$ is
a non-transitive ergodic flow and so is $F^M$.
Therefore, $M$ is of type $\mathrm{III}_0$.
\end{proof}

We give a simple example of a minimal action of a compact group with
non-trivial $m_\beta(\hK)$.
Let $\lambda$ be a positive number with $0<\lambda<1$, and set
$T=-2\pi/\log \lambda$.
We define a state of the 3 by 3 matrix algebra $M(3,\C)$ by
$$\varphi_0(x)=\mathrm{Tr}(\Di(\frac{1}{1+2\lambda},
\frac{\lambda}{1+2\lambda},\frac{\lambda}{1+2\lambda})x), $$
where $\Di$ means a diagonal matrix.
We set
$$(M,\varphi)=\bigotimes_{n=1}^\infty (M(3,\C),\varphi_0),$$
which is a Powers factor of type $\mathrm{III}_\lambda$.
For $g\in SU(2)$, we define
$$\alpha_g=\bigotimes_{n=1}^\infty (\Ad(\Di(1,g))).$$
Then, $\alpha$ is a minimal action of $SU(2)$.
It is easy to show
$\sigma^{\varphi}_{T/2}=\alpha_{-1}$, which implies
that $P=M^{\alpha_{-1}}$ is of type $\mathrm{III}_{\lambda^2}$.
Note that the restriction of $\alpha$ to $P$ induces a minimal
action of $SO(3)$, whose center is trivial.
Thus, Corollary 5.14 implies that $M^\alpha$ and $P$ has the common
flow and $M^\alpha$ is of type $\mathrm{III}_{\lambda^2}$.
Corollary 5.14 again implies that the endomorphism in the Roberts
action corresponding to the spin $1/2$ representation has a non-trivial
Connes-Takesaki module.

\appendix
\section{Miscellaneous results in ergodic theory}

In this appendix, we collect some results from ergodic theory,
(or its simple non-commutative generalization), used in the main body of
the present notes.
We provide proofs here because we could not find appropriate references,
though some (or all) of them are probably well-known to specialists.

For a locally compact group $G$, we say that $(X,\mu)$ is a $G$-space
if $X$ is a standard Borel space, $\mu$ is a probability measure on
$X$, and $G$ acts on $(X,\mu)$ as a non-singular Borel transformation group.
We say that a $G$-action is faithful if the corresponding action on
$L^\infty(X,\mu)$ is faithful.
We say that two $G$-spaces $(X_1,\mu_1)$ and $(X_2,\mu_2)$ are isomorphic
if the corresponding $G$-actions on the von Neumann algebras
(or equivalently, measure algebras)
$L^\infty(X_1,\mu_1)$ and $L^\infty(X_2,\mu_2)$ are conjugate.
For $g\in G$ and a Borel subset $E\subset X$, we use the notation
$g\cdot \mu(E):=\mu(E\cdot g)$.

\subsection{Compact automorphism groups of ergodic transformation groups}
\begin{lemma} Let $G$ be a locally compact group and
$(X,\mu)$ be an ergodic $G$-space.
We assume that a compact group $K$ faithfully acts on $(X,\mu)$ commuting
with $G$.
We use the notation with $G$ acting on $X$ from right and $K$ acting on $X$
from left.
Let $(Y,\nu)$ be the factor space of $(X,\mu)$ by $K$.
Then, there exists a minimal $K$-valued cocycle $c:Y\times G\longrightarrow
K$
with $K_c=K$ such that $(X,\mu)$ is isomorphic to the skew product
$(Y\times_c K,\nu\times m_K)$, that is,
under this isomorphism, the factor map from $X$ onto $Y$ corresponds to the
projection onto the first component, and the $K\times G$-action corresponds
to
the following action:
$$k\cdot (\omega,l)=(\omega,kl),\quad \omega\in Y,\; k,l\in K,$$
$$(\omega,l)\cdot g=(\omega\cdot g,lc(\omega,g)),\quad g\in G.$$
$c$ is unique up to equivalence.
\end{lemma}

\begin{proof} Considering the Gelfand spectrum of an appropriate weakly
dense
separable $C^*$-subalgebra of $L^\infty(X,\mu)$,  we may and do assume that
$X$ is a compact metric space, on which $G\times K$ acts continuously, and
that
$\mu$ is invariant under $K$.
Since $K$ is compact, every $K$-orbit is closed and there exists a Borel
subset
$Y\subset X$ that meets each $K$-orbit exactly once (see, for example,
\cite[Section 3.4]{A}).
Let $\cK$ be the set of all closed subgroups of $K$.
Then, $\cK$ has a Polish space structure such that $K$ continuously acts on
$\cK$ by $k\cdot H=kHk^{-1}$, $H\in \cK$, $k\in K$ \cite{Fe}.
We denote by $[H]$ the class of $H$ in the quotient space of $\cK/K$.
For each $\omega\in X$, let $K_\omega\subset K$ be the stabilizer subgroup
of $\omega$.
Then, thanks to \cite[Chapter II, Proposition 2.3]{AM}, the map
$$X\ni \omega\mapsto [K_\omega]\in \cK/K$$
is a Borel map satisfying $[K_{\omega\cdot g}]=[K_\omega]$ for $g\in G$.
Thanks to the ergodicity of the $G$-action, there exist $H\in \cK$ and
$K\times G$-invariant Borel null set $E\subset X$ such that for
every $\omega\in X_0:=X\setminus E$, we have $[K_\omega]=[H]$.
We set $Y_0=Y\cap X_0$.

Let $N(H)$ be the normalizer subgroup of $H$ in $K$.
Choosing a Borel cross section from $N(H) \backslash K$ to $K$, we have a
Borel map
$\varphi:Y_0\longrightarrow K$ such that
$\varphi(\omega)\cdot H=K_\omega$.
We set
$$Y_1=\{\varphi(\omega)^{-1}\cdot \omega\in X_0;\omega\in Y_0\}.$$
Then, $Y_1$ is a Borel subset of $X_0$ that meets each $K$-orbit exactly
once
and $K_\omega=H$ for all $\omega\in Y_1$.
Choosing a Borel cross section $\psi: H\backslash K\longrightarrow K$,
we set
$$\Phi(\omega,kH)=\psi(kH)\cdot \omega,\quad \omega \in Y_1,\; kH\in
H\backslash K.$$
Thanks to \cite[pp. 72, Corollary]{A}, $\Phi$ is a Borel isomorphism
from $Y_1\times H\backslash K$ onto $X_0$.
Since $\mu$ is $K$-invariant, $\mu\cdot \Phi$ is of the form
$\nu\times m_{H\backslash K}$, where $m_{H\backslash K}$ is the
(normalized) Haar measure of $H\backslash K$.
We introduce a $K\times G$-action on $Y_1 \times H\backslash K$ from the
action
on $X_0$ through $\Phi$.

When we regard $(Y_1,\nu)$ as the factor space of $X_0$ by the $K$-action,
we denote the action of $g\in G$ on $Y_1$ by $\omega \circ g$,
$\omega\in Y_1$.
Then, for each $\omega\in Y_1$ and $g\in G$, there exists a unique
$c_{(\omega,g)}H\in H\backslash K$ such that
$\omega\cdot g=c_{(\omega,g)}\cdot (\omega\circ g)$.
Using a Borel cross section from $H\backslash K$ to $H$ again, we can take
$c(\omega,g)\in K$ to be a Borel map from $Y_1\times G$ to $K$.
Let $\pi:K\longrightarrow H\backslash K$ be the quotient map.
Commutativity of the $K$-action and $G$-action implies that for each fixed
$g\in G$, $c$ takes its values in $N(H)$ almost everywhere.
Moreover, $\pi\cdot c$ is an $N(H)/H$-valued cocycle and
$$(\omega, \pi(k))\cdot g=(\omega\circ g, \pi(kc(\omega,g)))
,\quad \omega\in Y_1,\; g\in G,\; k\in K. $$
This implies that for every subset $U\subset K$,
$Y_1\times \pi(U\cdot N(H))$ is $G$-invariant.
Thus, thanks to the ergodicity of the $G$-action, we conclude that
whenever $U$ is an open subset of $K$,
$$m_K(U\cdot N(H))=m_{H\backslash K}(\pi(U\cdot N(H)))=1,$$
which implies $N(H)=K$, that is, $H$ is a normal subgroup of $K$.
Since the $K$-action is faithful, we get $H=\{e\}$.

Uniqueness of $c$ up to equivalence is obvious.
\end{proof}

\subsection{Non-commutative extensions with relatively discrete spectrum}
In \cite{Z1}, Zimmer introduced the notion of ergodic extensions with
relatively discrete spectrum.
The following is a slight generalization of this notion to
the non-commutative case.

\begin{definition} Let $G$ be locally compact group and $(X,\mu)$ be
an ergodic $G$-space.
Let $A$ be a von Neumann algebra with an ergodic $G$-action $\Phi$
such that the center of $A$ equivariantly includes $Z=L^\infty(X,\mu)$,
that is
$$\Phi_g(f)(\omega)=f(\omega\cdot g),\quad f\in Z,\; \omega\in X,\;
g\in G.$$
Let $u$ be the canonical implementation of $\Phi$ on the standard
Hilbert space for $A$.
We say that $A$ is an \textit{extension of $Z$ with relatively discrete
spectrum} if the von Neumann algebra $u(G)'\cap Z'$ is decomposed into
a direct sum of type I factors, and for each minimal projection
$p\in u(G)'\cap Z'$, $\natural(p)$ is bounded, where $\natural$ is
the (non-normalized) center valued trace on $Z'$, (in other words,
the restriction of $\natural$ to $u(G)'\cap Z'$ is semifinite).
\end{definition}

Note that the above definition coincides with Zimmer's one
\cite[Definition 4.2]{Z1} when $A$ is commutative
(although the two definitions appear different).

\begin{lemma} Let $G$ be a locally compact group, $K$ be a compact group,
$(X,\mu)$ be a $G$-space, and $c:X\times G \longrightarrow K$ be a cocycle.
If $\Psi$ is an action of $K$ on a von Neumann algebra $B$ standardly
acting on a Hilbert space $H$, then there exists an action $\Psi^c$ of
$G$ on $L^\infty(X,\mu)\otimes B$ such that for all $g\in G$,
$f\in L^\infty(X,\mu)$, $x\in B$, $\xi\in L^2((X,\mu),H)$,
$$(\Psi^c_g(f\otimes x)\xi)(\omega)
=f(\omega\cdot g)\Psi_{c(\omega,g)}(x)\xi(\omega)$$
holds for almost every $\omega\in X$.
\end{lemma}

\begin{proof} Let $v$ be the canonical implementation of $\Psi$ on $H$.
We introduce a unitary $u(g)$, $g\in G$ on  $L^2((X,\mu),H)$ by
$$(u(g)\xi)(\omega)=\frac{dg\cdot\mu}{d\mu}(\omega)^{1/2}
v(c(\omega,g))\xi(\omega\cdot g),\quad \xi\in L^2((X,\mu),H),\;
\omega\in X.$$
Then, $u$ is a representation of $G$, which is a Borel map
from $G$ to the unitary group of $B(L^2((X,\mu), H))$ with respect to
the Borel structure coming from the weak (and strong ) topology.
Since any unitary Borel representation of a locally compact group
is continuous, (or more strongly, any Borel homomorphism between Polish
groups
is continuous \cite[Chapter I, Propositon 3.3]{AM}), so is $u$.
$\Psi^c_g$ is defined by the restriction of $\Ad(u(g))$ to
$L^\infty(X,\mu)\otimes B$.
\end{proof}

The proof of (1) of the following theorem is a straightforward
generalization
of
Zimmer's argument in \cite[Theorem 4.3]{Z1} except for the proof of
uniqueness,
which does not exist there.
(Zimmer's comment about uniqueness in \cite[pp.401]{Z1} does not seem to be
relevant.)

\begin{theorem} Let $G$ be a locally compact group, $(X,\mu)$ be an ergodic
$G$-space, and $\Phi$ be an action of $G$ on a von Neumann algebra $A$.
We assume that  $A$ is an extension of $Z:=L^\infty(X,\mu)$ with relatively
discrete spectrum. Then,
\begin{itemize}
\item[$(1)$] There exist a compact group $K$, a faithful ergodic action
$\Psi$ of $K$ on a von Neumann algebra $B$, and a minimal cocycle
$c:X\times G \longrightarrow K$ with $K_c=K$ such that
$$(Z\subset A,\Phi)\cong
(Z\otimes \C\subset Z\otimes B,\Psi^c),$$
that is, there exists an isomorphism from $A$ onto $Z\otimes B$ that
intertwines $\Phi$ and $\Psi^c$, and is identity on $Z$.
\item[$(2)$] $A$ is always an injective von Neumann algebra, and there
exists
a unique faithful normal conditional expectation $\epsilon$ from $A$ onto
$Z$
satisfying $\Phi_g\cdot \epsilon= \epsilon\cdot\Phi_g$, $g\in G$ and
$\epsilon(ab)=\epsilon(ba)$, $a,b\in A$.
\item[$(3)$] $(K,c, \Psi, B)$ satisfying the condition of (1) are unique
in the following sense: if $(K^1,c^1, \Psi^1, B^1)$ satisfy the same
condition, then
there exist a group isomorphism $\nu$ from $K$ onto $K_1$ and an isomorphism
$\sigma$ from $B$ onto $B^1$  such that $\nu\cdot c$ and $c^1$ are
equivalent
and $\Psi^1_{\nu(k)}\cdot \sigma=\sigma\cdot \Psi_k$, $k\in K$.
\end{itemize}
\end{theorem}

\begin{proof} (1). By assumption, we have a direct sum decomposition
$$Z'\cap u(G)'=\bigoplus_{\lambda\in \Lambda}M_\lambda,$$
where $M_\lambda$  is a type I factor.
We take a system of matrix units $\{e^\lambda_{ij}\}_{i,j=1}^{n_\lambda}$ of
$M_\lambda$ with $n_\lambda$ either a natural number or possibly infinite
such that $\{e^\lambda_{ii}\}_{i=1}^{n_\lambda}$ form a partition of unity
consisting of minimal projections of $M_\lambda$.
Since $\Ad(u(g))$, $g\in G$ normalizes $Z$, it does $Z'$ as well, and we set
$\Phi'_g$ to be the restriction of $\Ad(u(g))$ to $Z'$.
Then, we have $\Phi_g\cdot \natural =\natural \cdot \Phi'_g,$
and in particular,
$$\Phi_g(\natural(e^\lambda_{11}))=\natural(e(\lambda)_{11}),
\quad \forall g\in G.$$
Thus, by ergodicity $\natural(e^\lambda_{11})$ is a constant, say
$d(\lambda)\in \N$.

In what follows, when we claim that ``some statement involving $g\in G$ and
$\omega\in X$ holds for almost all $\omega\in X$", we always mean that for
each fixed $g\in G$, there exists a conull set $E_g\subset X$ such that
the statement holds for all $\omega\in E_g$.
The annoying fact that $E_g$ varies according to $g\in G$ is usually
taken care of by Mackey's argument in \cite{M1}.
We consider the disintegration of $H$ over $Z$
$$H=\int_{X}^\oplus H(\omega)d\mu(\omega).$$
Then, $\natural(e^\lambda_{11})=d(\lambda)$ means that the dimension of
$e^\lambda_{11}(\omega)H(\omega)$ is $d(\lambda)$ for almost all $\omega\in
X$.
We take $\xi^\lambda_a=(\xi^\lambda_a(\omega))\in H$,
$a=1,2,\cdots,d(\lambda)$ such that for almost all
$\omega\in X$ $\{\xi^\lambda_a(\omega)\}_{a=1}^{d(\lambda)}$ is an
orthonormal
basis of $e^\lambda_{11}(\omega)H(\omega)$.
Since $u(g)$ is a unitary satisfying $(u(g)fu(g)^*)(\omega)=f(\omega\cdot
g)$
for $f\in Z$, for $\xi,\eta\in H$ and $f\in Z$ we have, on one hand,
$$\inpr{u(g)f\xi}{u(g)\eta}=\int_Xf(\omega\cdot g)
\inpr{(u(g)\xi)(\omega)}{(u(g)\eta)(\omega)}d\mu(\omega),$$
and on the other hand,
\begin{eqnarray*}\inpr{u(g)f\xi}{u(g)\eta}&=&
\inpr{f\xi}{\eta}=\int_X f(\omega)
\inpr{\xi(\omega)}{\eta(\omega)}d\mu(\omega)\\
&=&\int_X f(\omega\cdot g)\frac{dg\cdot\mu}{d\mu}(\omega)
\inpr{\xi(\omega\cdot g)}{\eta(\omega\cdot g)}d\mu(\omega).
\end{eqnarray*}
This implies,
$$\inpr{(u(g)\xi)(\omega)}{(u(g)\eta)(\omega)}
=\frac{dg\cdot\mu}{d\mu}(\omega)
\inpr{\xi(\omega\cdot g)}{\eta(\omega\cdot g)},$$
for almost all $\omega$.
Thus, $$\{\frac{dg\cdot\mu}{d\mu}(\omega)^{-1/2}
(u(g)\xi^\lambda_a)(\omega)\}_{i=1}^{d(\lambda)}$$
is an orthonormal basis of $e^\lambda_{11}(\omega)H(\omega)$ for almost all
$\omega$.
We set
$$c^\lambda(\omega,g)_{ab}=\frac{dg\cdot\mu}{d\mu}(\omega)^{-1/2}
\inpr{(u(g)\xi^\lambda_a)(\omega)}{\xi^\lambda_b(\omega\cdot g)}.$$
Then, $c^\lambda(\omega,g)=(c^\lambda(\omega,g)_{ab})$ is an element of
$Z^1(G,U(d(\lambda)))$.
We define a compact group $\cK$ and a cocycle $c$ by
$$\cK=\prod_{\lambda\in \Lambda}U(d(\lambda)),$$
$$c=(c^\lambda)_{\lambda\in \Lambda}\in Z^1(G,\cK).$$
Changing the system $\{\xi^\lambda_a\}_{a=1}^{d(\lambda)}$ if necessary,
we may and do assume that $c$ is minimal.
We denote by $K$ the minimal subgroup of $c$.
Note that $c^\lambda$ is minimal for all $\lambda\in \Lambda$ and that if
$\lambda_1\neq\lambda_2$, $c_{\lambda_1}$ and $c_{\lambda_2}$ are not
equivalent (otherwise, $Z'\cap u(G)'$ would be bigger).

Let $H_\lambda^1$ and $H_\lambda^2$ be Hilbert spaces with basis
$\{f^\lambda_a\}_{a=1}^{d(\lambda)}$ and $\{e^\lambda_i\}_{i=1}^{n_\lambda}$
respectively.
We set
$$H_0=\bigoplus_{\lambda\in \Lambda}H_\lambda^1\otimes H_\lambda^2.$$
Since $\{(e^{\lambda}_{i1}\xi^\lambda_a)(\omega)\}_{\lambda,i,a}$ is an
orthonormal basis of $H(\omega)$ for almost all $\omega$,
we identify $(e^{\lambda}_{i1}\xi^\lambda_a)(\omega)$ with
$f^\lambda_a\otimes e^\lambda_i$, and $H$ with $L^2(X,\mu)\otimes H_0$.
There exists a natural embedding of
$\cK=\prod_{\lambda\in \Lambda} U(H_\lambda^1)$ in the group $U(H_0)$ of
all unitaries $H_0$ as those unitaries acting ``only on $H_\lambda^1$,
$\lambda\in \Lambda$".
Through this embedding, we also regard $K$ as a subgroup of $U(H_0)$.
Then, we have
\begin{equation}
(u(g)\xi)(\omega)=(\frac{dg\cdot\mu}{d\mu}(\omega))^{1/2}
c(\omega,g)\xi(\omega\cdot g),\quad \xi\in L^2((x,\mu),H_0),
\end{equation}
for almost all $\omega\in X$.
Let
$$A=\int_X^{\oplus}A(\omega)d\mu(\omega)$$
be the disintegration of $A$ over
$Z=L^\infty(X,\mu)\otimes \C1_{H_0}$.
Then, the map $X\ni\omega\mapsto A(\omega)$ is measurable
with respect to the Effros Borel structure of the set of von Neumann
algebras on $H_0$ \cite[Chapter 4]{N},
\cite[Chapter 4, Section 8]{T2}.
For $x\in A$, (A.1) implies
\begin{equation}
\alpha_g(x)(\omega)=c(\omega,g)x(\omega\cdot g)c(\omega,g)^{-1},
\end{equation}
for almost all $\omega$.
Since $A$ has countable generators (we assume that Hilbert spaces are
separable throughout this paper), this implies
\begin{equation}
c(\omega,g)A(\omega\cdot g)c(\omega,g)^{-1}=A(\omega),
\end{equation}
for almost all $\omega$.

Using the Effros Borel structure in the same way as in the proof of
\cite[Theorem 4.3]{Z1} and passing to an equivalent cocycle if necessary
(and changing $A(\omega)$ and $c(\omega,g)$ on a null set),
we may and do assume that there exists a von Neumann algebra
$B\subset B(H_0)$ such that $A(\omega)=B$ and $c(\omega,g)$ normalizes
$B$ for all $\omega\in X$.
Since $c$ is a minimal cocycle with $K_c=K$, $K$ normalizes $B$, and so
we define $\Psi$ to be the restriction of $\Ad(k)$, $k\in K$ to $B$.
Then, (A.2) implies $\Phi=\Psi^c$.
Since $\Phi$ is ergodic, so is $\Psi$.

(2). Thanks to \cite{HLS}, $B$ is injective and so is $A$.

Let $\epsilon$ be a normal conditional expectation from $A$ onto $Z$.
Then, since $H$ is the standard Hilbert space of $A$, there exists
$\xi=(\xi(\omega))\in L^2((X,\mu),H_0)$ whose vector state is
$\mu\cdot \epsilon$, where $\mu$ means the state of $Z$ corresponding
to the measure $\mu$.
Thus, for all $x\in A$ and $f \in Z$ we have
\begin{equation}
\int_Xf(\omega)\epsilon(x)(\omega)d\mu(\omega)=\inpr{\epsilon(fx)\xi}{\xi}.
\end{equation}
We introduce a measurable field of normal states $\{\varphi_\omega\}$ of $B$
by
$$\varphi_\omega(b)=\inpr{b\xi(\omega)}{\xi(\omega)},\quad b\in B.$$
Then, (A.4) implies $\epsilon(x)(\omega)=\varphi_{\omega}(x(\omega))$
for almost all $\omega$.
Moreover, if $\epsilon$ commutes with $\Phi_g$, we get
$$\varphi_{\omega\cdot g}=\varphi_\omega\cdot\Psi_{c(\omega,g)},$$
for almost all $\omega$.
Let $\mathfrak{S}(B)$ be the set of all normal states of $B$,
which is a Polish space in the norm topology.
Since $\Psi$ is continuous in $u$-topology, $K$ acts on $\mathfrak{S}(B)$
as a continuous transformation group.
Applying the same argument as above to $\mathfrak{S}(B)$ and passing
to an equivalent cocycle (and changing $\varphi_\omega$ and
$c(\omega,g)$ on a null set), we may assume that there exists
$\varphi \in \mathfrak{S}(B)$ such that $\varphi_\omega=\varphi$ for
all $\omega$ and $\varphi\cdot \Psi_{c(\omega,g)}=\varphi$ holds.
Thus, $\varphi$ is a $\Psi$-invariant state, which is unique because of
ergodicity of $\Psi$.
Thus, $\epsilon=\id\otimes \varphi$ is uniquely determined.
$\epsilon$ has the trace property as $\varphi$ is a trace thanks to
\cite{HLS}.

(3). We assume that $(A, \Phi)$ are also given by
$A=Z\otimes B^1$, $\Phi=(\Psi^1)^{c^1}$.
Let $\tau$ be the unique $\Psi^1$-invariant trace, $L^2(B^1)$
be the GNS Hilbert space of $\tau$, and $v(k)$ be the canonical
implementation of $\Psi^1_k$, $k\in K^1$.
Then, we have $H=L^2(X,\mu)\otimes L^2(B^1)$ and $u$ is given by
$$u(g)\xi(\omega)=\frac{dg\cdot \mu}{d\mu}(\omega)^{1/2}v(c^1(\omega,g))
\xi(\omega\cdot g),\quad \xi\in L^2((X,\mu),L^2(B^1)).$$
Since $K^1$ is a compact group, we have the irreducible decomposition
$$(L^2(B^1),v(k))=(\bigoplus_{(\pi,H_\pi)\in \hK^1} H_\pi\otimes L_\pi,
\bigoplus_{(\pi,H_\pi)\in \hK^1} \pi(k)\otimes 1_{L_\pi}),$$
where $L_\pi$ is the multiplicity space, which may be zero.
Since the $G$-action on $(X,\mu)$ is ergodic and $c^1$ is a minimal cocycle
with
$K_{c^1}=K$ we have
$$Z'\cap u(G)'=\bigoplus_{(\pi,H_\pi)\in \hK^1} \C\otimes \C\otimes
B(L_\pi).$$
Thus, we identify $\Lambda$ with the subset of $\hK^1$ consisting of the
irreducibles appearing in $(v,L^2(B^1))$, and identify $H^2_\lambda$ with
$L_\pi$ when $\lambda\in \Lambda$ is identified with $(\pi,H_\pi)\in \hK^1$.
In what follows, we abuse notation and use only $\Lambda$, instead of
$\hK^1$, if there is no possibility of confusion.
We set
$$\cK_1=\prod_{\lambda\in \Lambda}U(H_\lambda),$$
and regard it as a subgroup of the unitary group $U(H)$ as before.
Since $\Psi^1$ is faithful, we identify $K_1$ with $v(K_1)\subset \cK_1$.

Let $w(\omega):H_0\longrightarrow L^2(B^1)$ be the measurable field of
unitaries
describing the identity map on $H$ with respect to the two
distinct splitting $H=L^2(X,\mu)\otimes H_0$ and $H=L^2(X,\mu)\otimes
L^2(B^1)$.
$w(\omega)$ is decomposed as
$$w(\omega)=\bigoplus_{\lambda\in \Lambda}w_\lambda(\omega)
\otimes 1_{H^2_\lambda},$$
where $w_\lambda(\omega):H^1_\lambda\longrightarrow H_\lambda$ is a
measurable
field of unitaries.
Then, we have
$$B^1=w(\omega)Bw(\omega)^*,$$
$$c^{1}(\omega,g)=w(\omega)c(\omega,g)w(\omega\cdot g)^*,$$
for almost all $\omega\in X$.
Thanks to uniqueness of the minimal subgroup up to conjugacy,
there exists a unitary $w^0_\lambda$ from $H^1_\lambda$ to $H_\lambda$ such
that
if we set
$$\nu_0=(\prod_{\lambda\in \Lambda}\Ad(w^0_\lambda)):\cK\longrightarrow
\cK_1,$$
then $\nu_0(K)=K_1$.
We denote by $\nu_1$ the restriction of $\nu_0$ to $K$.
Let $w^1(\omega)=w_\lambda(\omega){w^0_\lambda}^*$, and
$$w^0=\bigoplus_{\lambda\in \Lambda} w^0_\lambda \otimes 1_{H^2_\lambda},$$
$$w^1(\omega)=\bigoplus_{\lambda\in \Lambda}w^1_\lambda(\omega)
\otimes 1_{H^2_\lambda},$$
Then, we have
$$B^1=w^1(\omega)w^0B{w^0}^*w^1(\omega)^*,$$
$$c^{1}(\omega,g)=w^1(\omega) \nu^1\cdot c(\omega,g) w^1(\omega\cdot g)^*,$$
for almost all $\omega\in X$.
Since $c^1$ and $\nu_1\cdot c$ are mutually equivalent minimal
$\cK_1$-valued
cocycles with $(\cK_1)_{c^1}=(\cK_1)_{\nu\cdot c}=K_1$,
\cite[Theorem 6.1]{Z1} implies that there exist $\kappa \in \cK_1$
normalizing
$K_1$ and a Borel map $a:X \longrightarrow K_1$ such that
$w^1(\omega)=a(\omega)\kappa$ for almost all $\omega\in X$.
Let $\nu(k)=\kappa \nu_1(k)\kappa^{-1}$, then we get
$$B^1=\kappa w^0B{w^0}^*\kappa^*,$$
$$c^{1}(\omega,g)=a(\omega)\nu\cdot c(\omega,g) a(\omega\cdot g)^*,$$
for almost all $\omega$.
Let $\sigma$ be the restriction of $\Ad(\kappa w^0)$ to $B$.
Then $\nu$ and $\sigma$ have the desired property.
\end{proof}

\subsection{Existence of minimal cocycles on ergodic flows}
The following is an easy generalization of \cite[Theroem 2]{Z2} to the
non-singular case, which asserts that there are plenty of non-trivial
modular endomorphisms for any type $\mathrm{III}_0$ factors.

\begin{proposition} Let $(X,\mu)$ be an ergodic non-transitive $\R$-space.
Then, for every compact group $K$, there exists a minimal cocycle
$c: X\times \R\longrightarrow K$ with $K_c=K$.
\end{proposition}

\begin{proof} As in \cite[Theorem 2]{Z2}, it suffices to show the statement
for a $\Z$-action instead of the $\R$-action thanks to Ambrose-Kakutani's
theorem \cite{AK}.
Let $S$ be a given non-singular and non-transitive ergodic transformation on
$(X,\mu)$.
We show that there exists a minimal cocycle
$c:X\times \Z\longrightarrow K$ with $K_{c}=K$.
We take a non-transitive ergodic transformation $T$ on $(Y,\nu)$
such that $\nu$ is $T$-invariant.
Thanks to \cite[Theorem 1]{Z2}, there exists a minimal cocycle
$c':Y \times \Z \longrightarrow K$ such that $K_{c'}=K$.
Let $R$ be the $\Z^2$-action on $X\times Y$ given by
$R_{(n_1,n_2)}=S^{n_1}\times T^{n_2}$.
We define a cocycle $c''$ on $(X\times Y)\times \Z^2$ by
$$c''((\omega_1,\omega_2),(n_1,n_2))=c'(\omega_2,n_2),$$
which is a minimal cocycle with $K_{c''}=K$.
Since $(X,\mu,S)$ is orbit equivalent to
$(X\times Y,\mu\times \nu, R)$ \cite{Kr}, we get the result.
\end{proof}

\end{document}